\documentclass{amsart}
\usepackage{graphicx}

\begin{document}

\begin{center}
\huge Geometry and Design of Equiangular Spirals
\end{center}
\
\begin{flushright}
\Large
Kostantinos Myrianthis
\normalsize \\
\vskip 1 mm
\verb+58, Zan Moreas str.,+\\
\verb+ Athens, P.C.15231,Greece+\\
\verb+myrian@ath.forthnet.gr+ \vskip 3 mm
\end{flushright}

                                                                                      \textbf{ABSTRACT}\\

In an equiangular spiral, "the whorls continually increase in breadth and do so in a steady and unchanging ratio... It follows that the sectors cut off by successive radii, at equal vectorial angles, are similar to one another in every respect and that the figure may be conceived as growing continuously without ever changing its shape the while" as stated by Sir D'Arcy W. Thompson and quoted in [1, p.125]. I was fascinated since my early years with the shape of spirals and all their versions in nature. The mathematical modeling of them became a very attractive topic of study and research for me and more specifically, the geometrical conditions under which any quadrangle or triangle can be fitted into similar copies of itself and form an equiangular spiral. This formation gives the impression of a digital form of spiral, where every digit is a triangle or quadrangle following similarity laws, which can allow a multiplicity of design capabilities. The essence of this work appears in the present article and is related with the geometry and the design characteristics of equiangular spirals.

\newpage

\section{\textbf{INTRODUCTION}}

\begin{figure}[bht]
\includegraphics{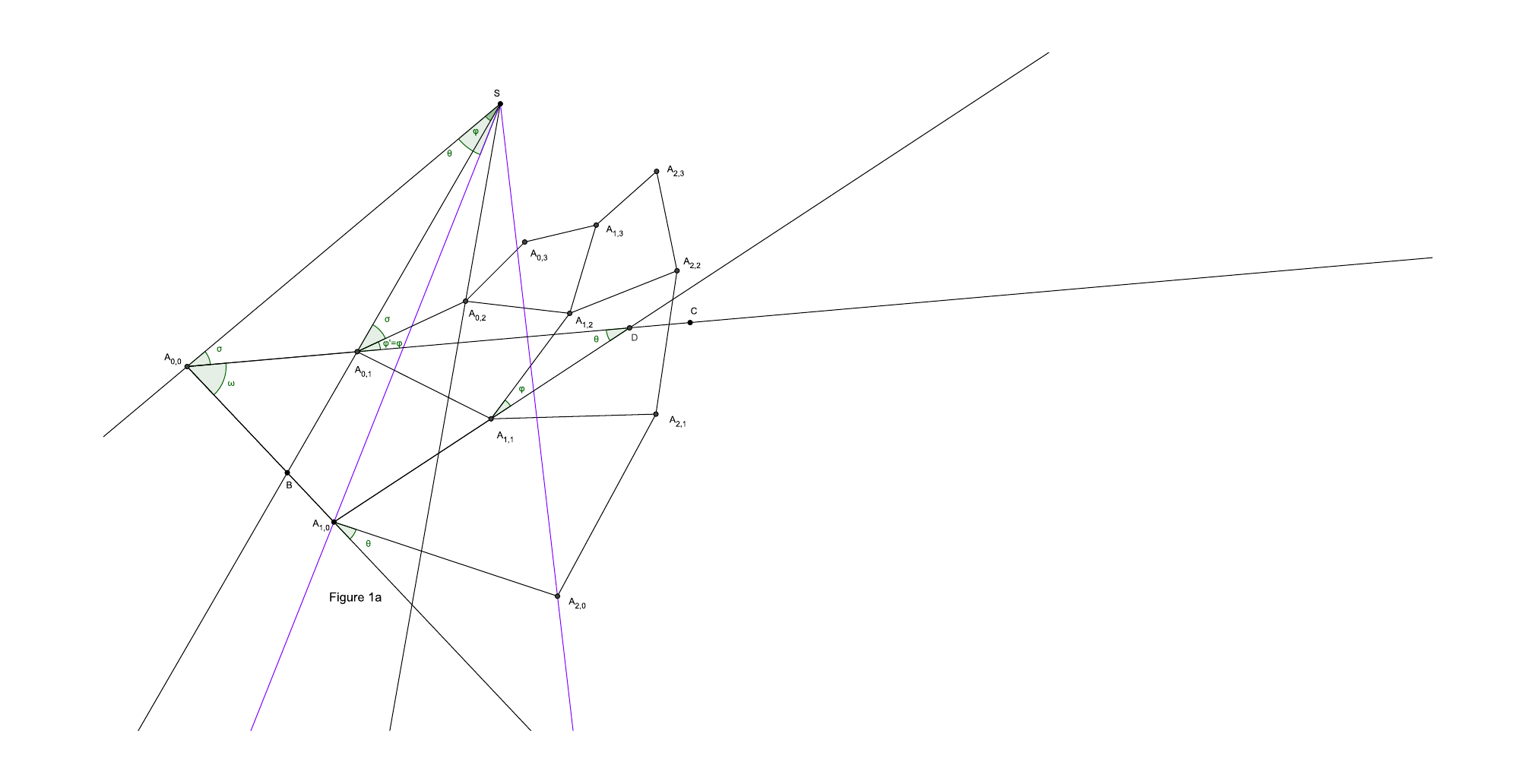}{}
\end{figure}

In figure 1a we have the branch $A_{0,0}A_{0,1}A_{0,2}A_{0,3}$..., defined by the equally spaced (by angle $\widehat{\phi}$) rays $SA_{0,0}$, $SA_{0,1}$, $SA_{0,2}$... from the centre $S$ of the spiral. The triangles $\triangle$$A_{0,0}A_{0,1}S$, $\triangle$$A_{0,1}A_{0,2}S$, $\triangle$$A_{0,2}A_{0,3}S$,... are similar with a ratio of similarity $\kappa$. As shown in Appendix 1, part1, $\widehat{(A_{0,0}A_{0,1},A_{0,1}A_{0,2})}=\widehat{\phi}$. The angle $\widehat{\sigma}$ can be expressed in relation to $\widehat{\phi}$ and $\kappa$ (as shown in Appendix 1, part 2):

\begin{equation}
 tan \widehat{\sigma}=\kappa sin \widehat{\phi} /(1-\kappa cos \widehat{\phi})
\end{equation}

\textbf{RULE 1}. It is obvious that given the position and size of 2 consecutive segments of a spiral branch so that $\widehat{\phi}$, $\widehat{\sigma}$  and $\kappa$ can be calculated, then the position of $S$ can be determined.\\

All the branches $A_{0,0}A_{0,1}A_{0,2}A_{0,3}$..., $A_{1,0}A_{1,1}A_{1,2}A_{1,3}$..., $A_{2,0}A_{2,1}A_{2,2}A_{2,3}$, can be named $B\kappa_{0}=A_{0,0}A_{0,1}A_{0,2}A_{0,3}...$, $B\kappa_{1}=A_{1,0}A_{1,1}A_{1,2}A_{1,3}...$,\\
$B\kappa_{2}=A_{2,0}A_{2,1}A_{2,2}A_{2,3}...$ or generally $B\kappa_{i}$ where $i=0,1,2,...$ in the direction of rotating of similarity ratio $\kappa$ for simplicity.\\

In the same figure 1a we have the branch $A_{0,0}A_{1,0}A_{1,1}A_{1,2}$..., defined by the equally spaced (by angle $\widehat{\theta}$) rays $SA_{0,0}$, $SA_{1,0}$, $SA_{2,0}$... from the centre $S$ of the spiral. The triangles $\triangle$$A_{0,0}A_{1,0}S$, $\triangle$$A_{1,0}A_{2,0}S$, $\triangle$$A_{2,0}A_{3,0}S$,... are similar with a ratio of similarity $\lambda$. As with the case of similarity ratio $\kappa$, $\widehat{(A_{0,0}A_{1,0},A_{1,0}A_{2,0})}=\widehat{\theta}$ and:

\begin{equation}
 tan (\widehat{\omega}+\widehat{\sigma})=\lambda sin \widehat{\theta} /(1-\lambda cos \widehat{\theta})
\end{equation}

As with the case of similarity ratio $\kappa$ all the branches $A_{0,0}A_{1,0}A_{2,0}A_{3,0}$...,\\
$A_{0,1}A_{1,1}A_{2,1}A_{3,1}$..., $A_{0,2}A_{1,2}A_{2,2}A_{3,2}$, can be named $B\lambda_{i}$ where $i=0,1,2,...$ in the direction of rotating of similarity ratio $\lambda$ for simplicity.\\

\textbf{RULE 2}. If the angle $\widehat{\sigma}$ is less than $\pi/2-\widehat{\phi}/2$, ($SA_{0,0}>SA_{0,1}$) then the branches $B\kappa_{i}$ rotate anticlockwise and converge to $S$, otherwise ($SA_{0,0}<SA_{0,1}$) they rotate clockwise and converge to $S$ (for $\widehat{\sigma}=\pi/2-\widehat{\phi}/2$ the triangle $\triangle$$A_{0,0}A_{0,1}S$ is isosceles and these branches do not converge). Similarly, if the angle $(\widehat{\omega}+\widehat{\sigma})$ is less than $\pi/2-\widehat{\theta}/2$ ($SA_{0,0}<SA_{1,0}$) then the branches $B\lambda_{i}$ rotate anticlockwise and converge to $S$, otherwise ($SA_{0,0}>SA_{1,0}$) they rotate clockwise and converge to $S$ (for $(\widehat{\omega}+\widehat{\sigma})=\pi/2-\widehat{\theta}/2$, these branches do not converge).\\

\begin{figure}[bht]
\includegraphics{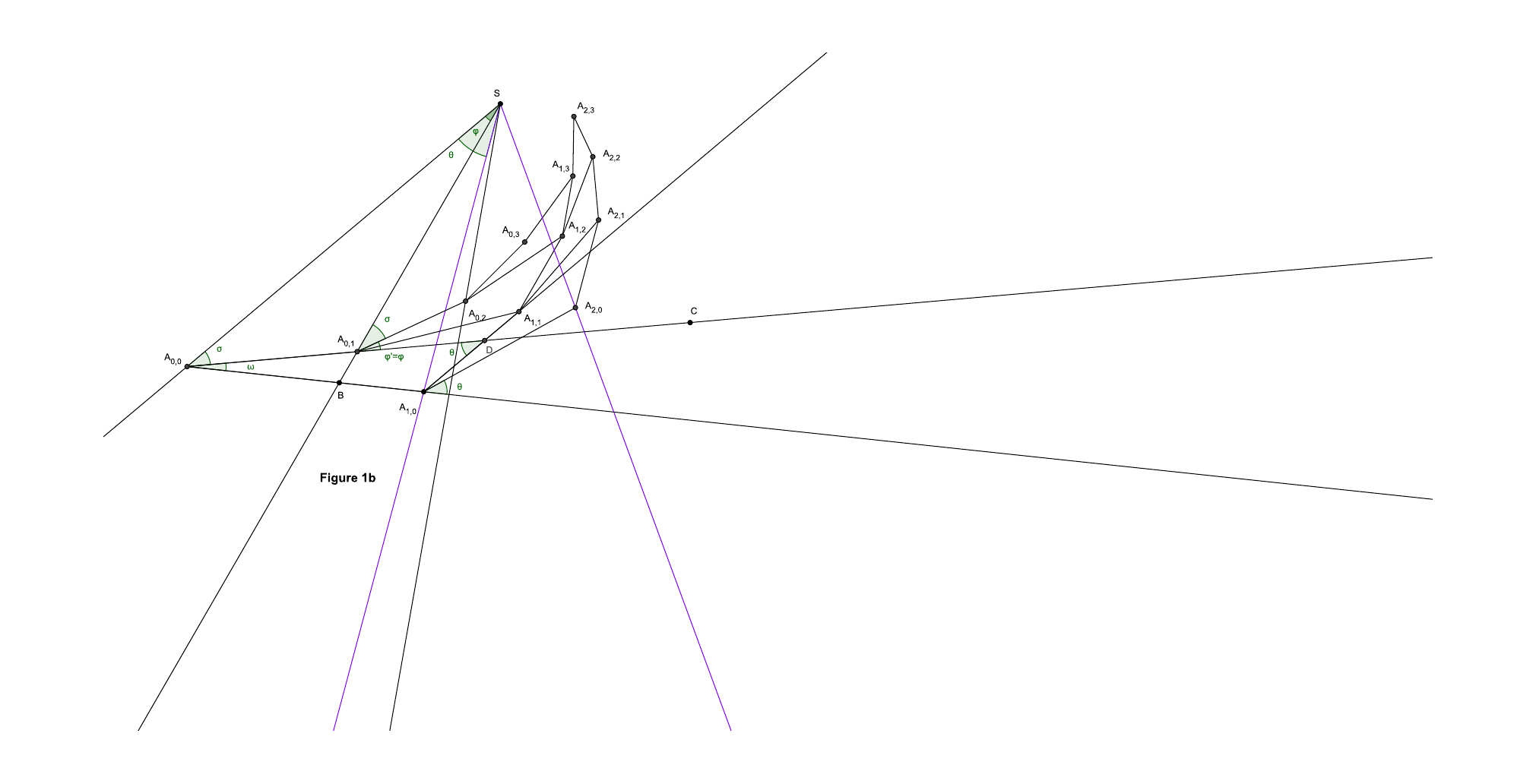}{}
\end{figure}

The quadrangles which are formed by the branches $B\kappa_{i}$ and $B_\lambda{i}$ (such as the initial one
$\Box A_{0,0}A_{0,1}A_{1,1}A_{1,0}$) can be either convex (figure 1a) or concave (figure 1b) and their relevant angles are calculated in Appendix 1, part 4. These quadrangles can become triangles in two cases, first when $\widehat{\omega}=\widehat{\phi}$ (figure 1c) and second when $\widehat{\omega}=0$ (figure 1d) and their relevant angles are also calculated in Appendix 1, part 4. It is interesting to note that for the case of $\widehat{\omega}=0$ of figure 1d, the triangles $\triangle$$A_{0,0}A_{1,0}S$, $\triangle$$A_{1,0}A_{2,0}S$,... which form the ratio $\lambda$ do not have a visual clear role. As shown in Appendix 1 part 5, we have the ratio $\overline{\lambda}=A_{0,2}A_{1,1}/A_{1,1}A_{2,0}=\kappa/\lambda$  related to the case of figure 1d, the triangles $\triangle$$A_{0,2}A_{1,1}S$, $\triangle$$A_{1,1}A_{2,0}S$ and the angle $\widehat{\overline{\theta}}=\widehat{A_{0,2}SA_{1,1}}=\widehat{A_{1,1}SA_{2,0}}=\widehat{\theta}-\widehat{\phi}$, where $\widehat{\theta}>\widehat{\phi}$ and these parameters do have a visual role (the case where $\widehat{\theta}<\widehat{\phi}$ can be treated accordingly). The parameters  $\widehat{\overline{\theta}}$ and $\overline{\lambda}$ can be treated just like $\widehat{\theta}$ and $\lambda$ at figure 1c.\\

\begin{figure}[bht]
\includegraphics{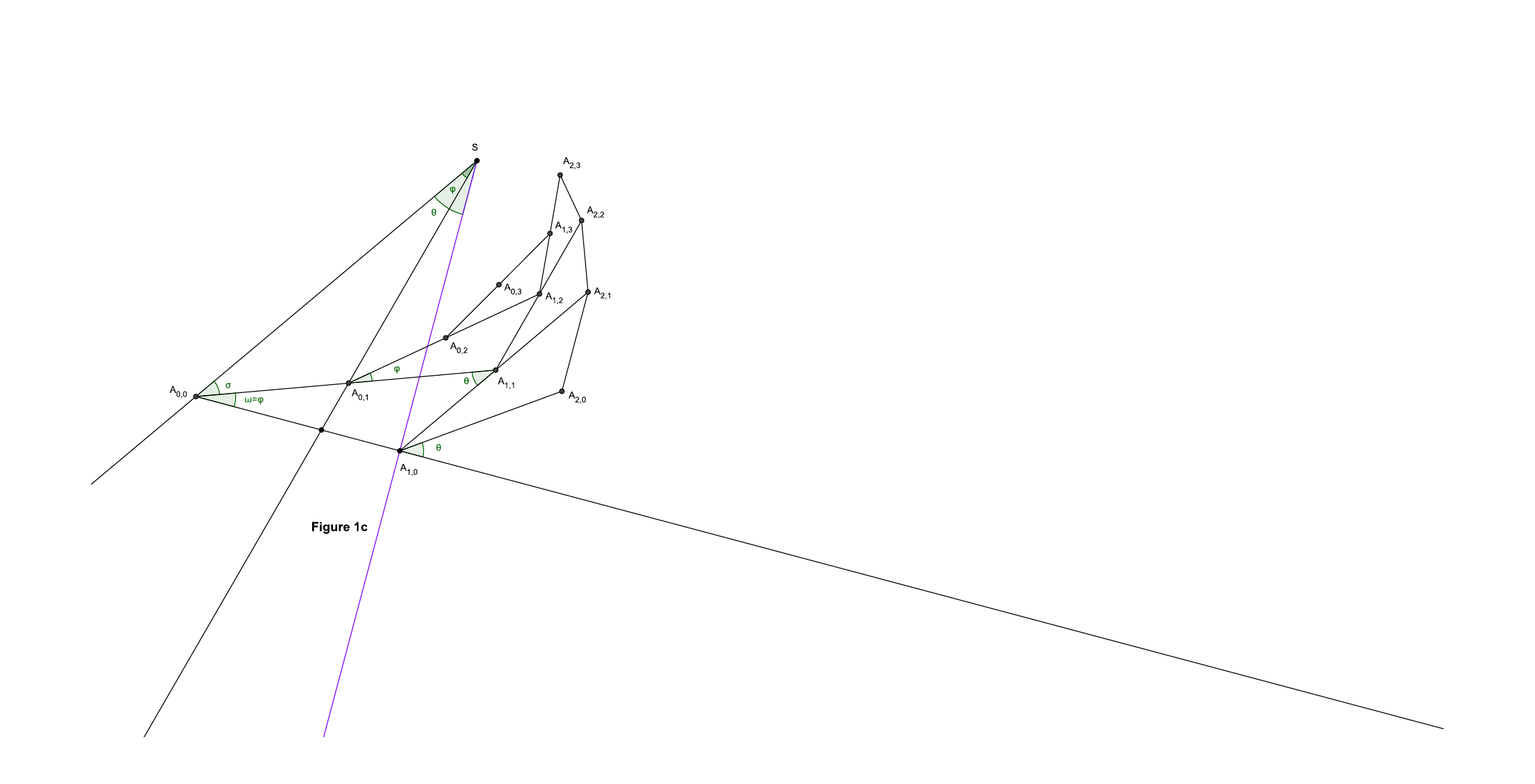}{}
\end{figure}

Starting from any quadrangle, for example the initial one, any quadrangle found in the $\kappa$ direction has a ratio of similarity $\kappa$ with the two adjacent quadrangles found in the same direction, therefore the quadrangles $\Box A_{0,0}A_{0,1}A_{1,1}A_{1,0}$ and $\Box A_{0,1}A_{0,2}A_{1,2}A_{1,1}$ have a ratio of similarity $\kappa$. Similarly, any quadrangle found in the $\lambda$ direction has a ratio of similarity $\lambda$ with the two adjacent quadrangles found in the same direction, therefore the quadrangles $\Box A_{0,0}A_{0,1}A_{1,1}A_{1,0}$ and $\Box A_{1,0}A_{1,1}A_{2,1}A_{2,0}$ have a ratio of similarity $\lambda$.\\

Also in figure 1a the branches $B\kappa_{i}$ are anticlockwise and $B\lambda_{i}$ are clockwise, so the branches in the $\kappa$ direction are contra-rotating with respect to the branches of the $\lambda$ direction. In the figures 1b, 1c, 1d, the branches $B\kappa_{i}$ and $B\lambda_{i}$ are co-rotating.\\

\begin{figure}[bht]
\includegraphics{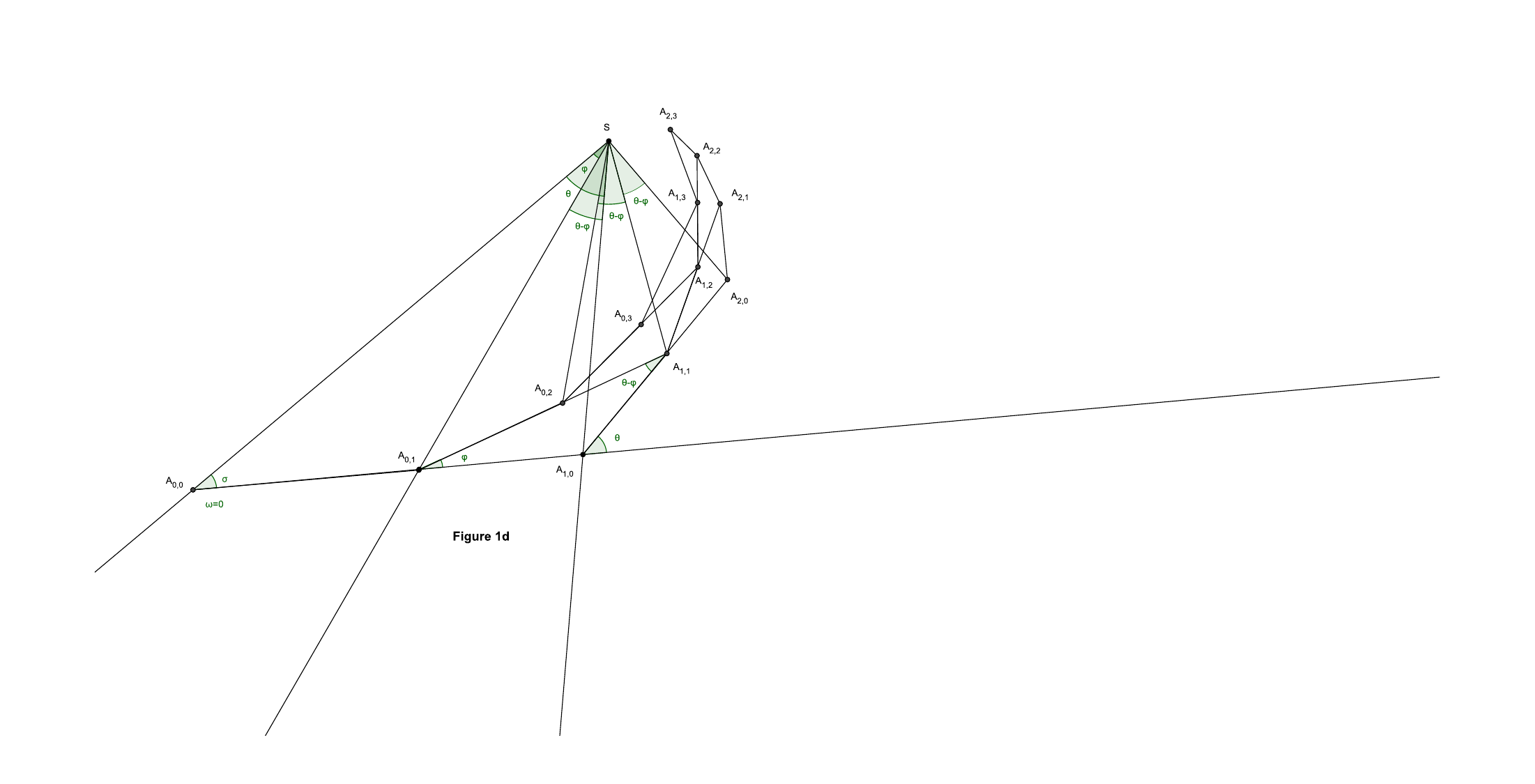}{}
\end{figure}

The basic parameters which define the geometry of a spiral system (with quadrangles or triangles) are $\widehat{\omega}$, $\widehat{\phi}$, $\widehat{\theta}$, $\kappa$ and $\lambda$. Given any four of them, the fifth can be calculated. Given these parameters and a segment of a spiral branch of the system which will be assumed as the starting segment such as $A_{0,0}A_{0,1}$ or $A_{0,0}A_{1,0}$, the spiral system can be constructed having a common starting ray $SA_{0,0}$ for both the first branch $B\kappa_{0}$ of the $\kappa$ direction and the first branch $B\lambda_{0}$ of the $\lambda$ direction.\\

Any arbitrary quadrangle (or triangle) which can be fitted into similar copies of itself, can produce a spiral system with branches of 2 types (one related with angle $\widehat{\phi}$ and ratio $\kappa$ and another one related with angle $\widehat{\theta}$ and ratio $\lambda$, as in figure 1a. The result usually overlaps itself or has cracks that cannot be filled, taking in account the related 2-D plane of the spiral system. These types of spiral systems (which can be called open), have an infinite number of branches in both $\kappa$ and $\lambda$ directions.\\

The intention of the present work is to provide the necessary conditions and parameter values using Euclidean geometry and algebra in order to produce a planar closed spiral system without cracks and overlaps, a closed spiral tiling system if one considers that quadrangles (or triangles) as tiles. In section 2, closed quadrangular systems are examined, whereas in sections 3 and 4 the two types of closed triangular systems ($\widehat{\omega}=\widehat{\phi}$ and $\widehat{\omega}=0$) are examined. All the possible cases of the equivalence  of triangular spirals are analyzed in section 5. Finally, in section 6, the concept of divergence angle ([4,p20]) is analyzed.\\

It is interesting to note that the topology of spiral systems is closely related to phyllotaxis theory as shown in [3] and [4] where complex exponential functions produce quadrangular and triangular closed spiral systems, which find wonderful applications in nature. A pure Euclidean approach such as the one of the present work might give a simpler and more practical point of view with enough interest for future development.

\section{\textbf{CLOSED QUADRANGULAR SPIRAL SYSTEMS}}

Let $\Box A_{0,0}A_{0,1}A_{1,1}A_{1,0}$ be one of the quadrangles shown in figure 2a, which according to the numbering system of the vertices of the quadrangles (or spiral branches) is the first one. Two branches appear in this figure, $B\lambda_{0}$ and $B\lambda_{1}$ in the $\lambda$ direction with 14 vertices each (0,..,13), and two branches $B\kappa_{0}$ and $B\kappa_{1}$ in the $\kappa$ direction with 4 vertices each (0,..,3), which co-rotate, (plus 12 branches $B\kappa_{i}$ where $i=2,...,13$ in the $\kappa$ direction with 2 vertices each 0,..,1). As it is shown in Appendix 1, part 3, $\widehat{(A_{0,0}A_{0,1},A_{1,0}A_{1,1})}=\widehat{\theta}$, therefore $\widehat{(A_{0,0}A_{0,1},A_{13,0}A_{13,1})}=13\widehat{\theta}$. Also it is obvious that $\widehat{(A_{0,0}A_{0,1},A_{0,2}A_{0,3})}=2\widehat{\phi}$.\\

\begin{figure}[bht]
\includegraphics{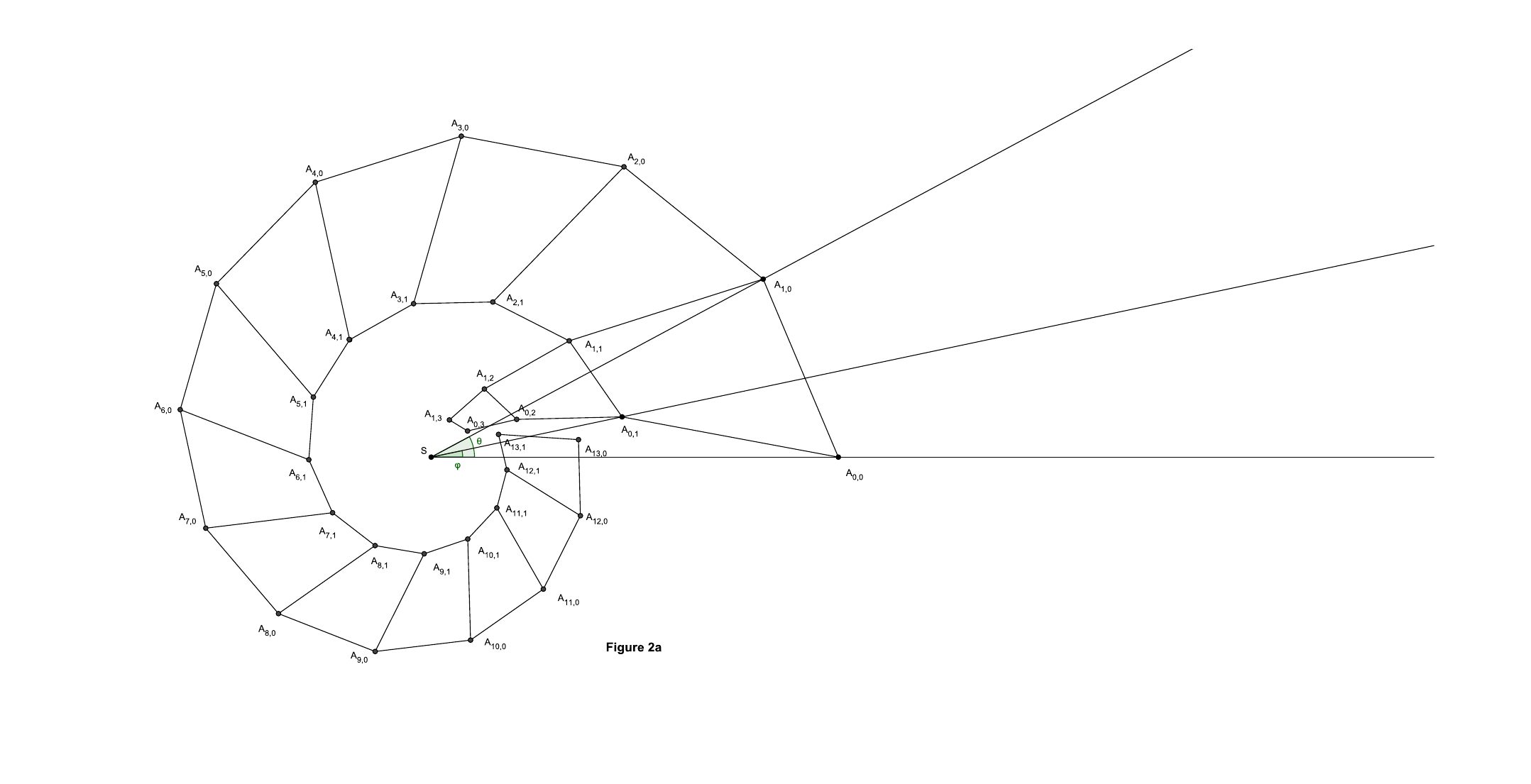}{}
\end{figure}

Obviously the open spiral system presented in figure 2a has a gap or it has overlaps if all the quadrangles are extended in the $\kappa$ and $\lambda$ directions. The initial branches $B\kappa_{0}$ and $B\lambda_{0}$ start from the same vertex $A_{0,0}$ and they never meet again at any vertex of the spiral, also there is an infinite number of branches in both $\kappa$ and $\lambda$ directions.\\

In figure 2b we have a spiral system of which the branches of $\kappa$ and $\lambda$ directions co-rotate, where there is no gap or overlap, so it can be considered as a closed spiral system. Because of this property, the numbering process of  the vertices (or nodes) of the quadrangles of the spiral system has to be different from the one of figure 2a. If we consider the segment $A_{13,1}A_{13,0}$ of figure 2a and assume that the spiral system of figure 2a is somehow transformed into the spiral system of figure 2b, it is as if this segment has moved and coincided with the segment $A_{0,2}A_{0,1}$ of figure 2a in order to produce the segment $A_{26,2}A_{13,1}$ of figure 2b. In other words, in figure 2b the initial branch $B\lambda_{0}$ starts from the vertex $A_{0,0}$ and meets for the first time the initial branch$B\kappa_{0}$ (which also starts from the vertex $A_{0,0}$) at the vertex $A_{13,1}$, which is the first vertex or point after the common start of the two branches. Given that a typical vertex of the spiral is $A_{i,j}$, in this case we have $i=13$ and $j=1$. Similarly, the two branches meet for the second time at the vertex $A_{26,2}$, and so on. In figure 2b we have one branch $B\lambda_{0}$ in the $\lambda$ direction and 13 branches $B\kappa_{i}$ where $i=0,..,12$ in the $\kappa$ direction.\\

\begin{figure}[bht]
\includegraphics{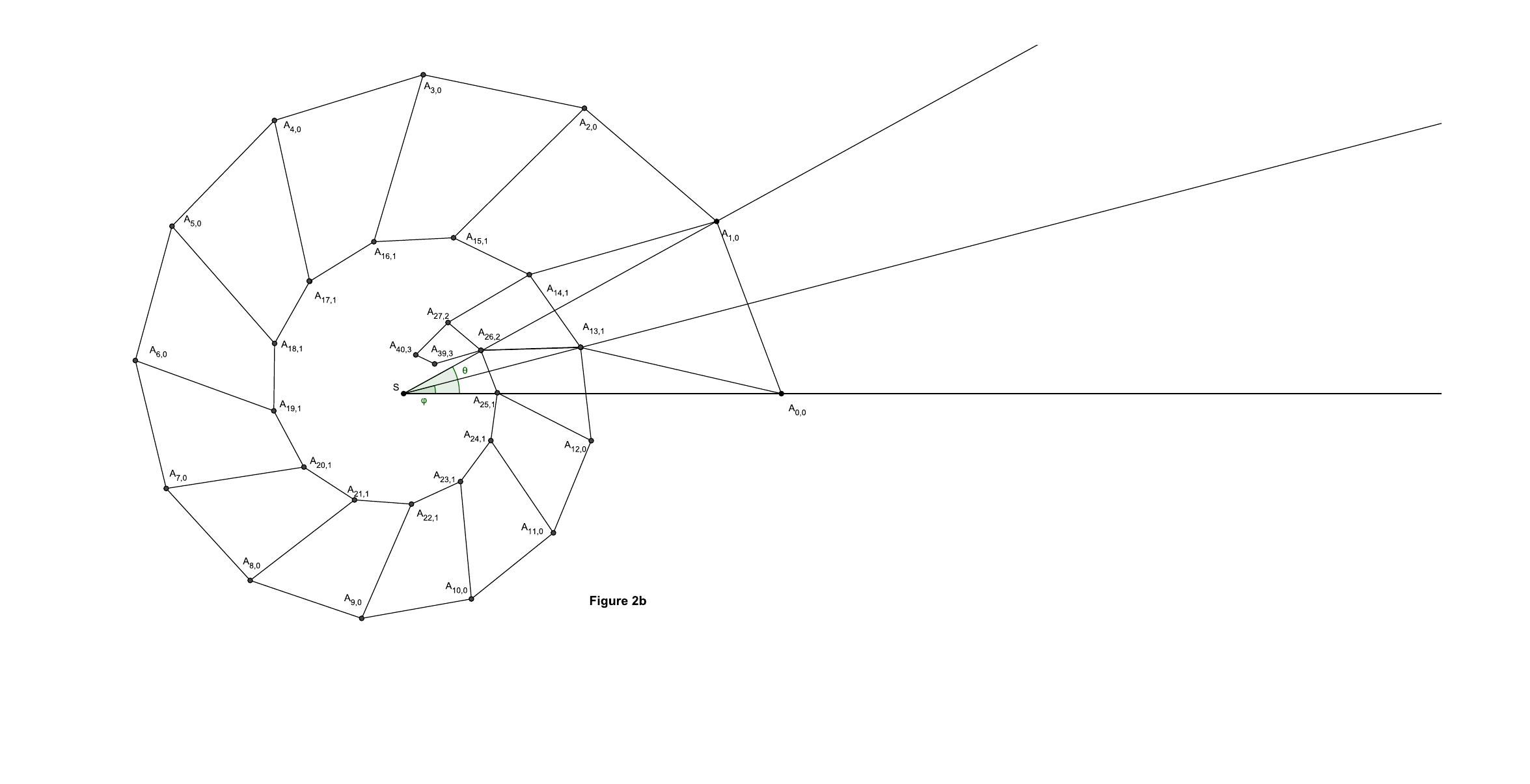}{}
\end{figure}

We have $\widehat{(A_{0,0}A_{13,1},A_{13,1}A_{26,2})}=\widehat{\phi}$ and $\widehat{(A_{0,0}A_{13,1},A_{1,0}A_{14,1})}=\widehat{\theta}$ (Appendix 1 part1 and part3). From the above and from figure 2b, we get:

\begin{equation}
 13\widehat{\theta}=2\pi+\widehat{\phi}
\end{equation}

\begin{figure}[bht]
\includegraphics{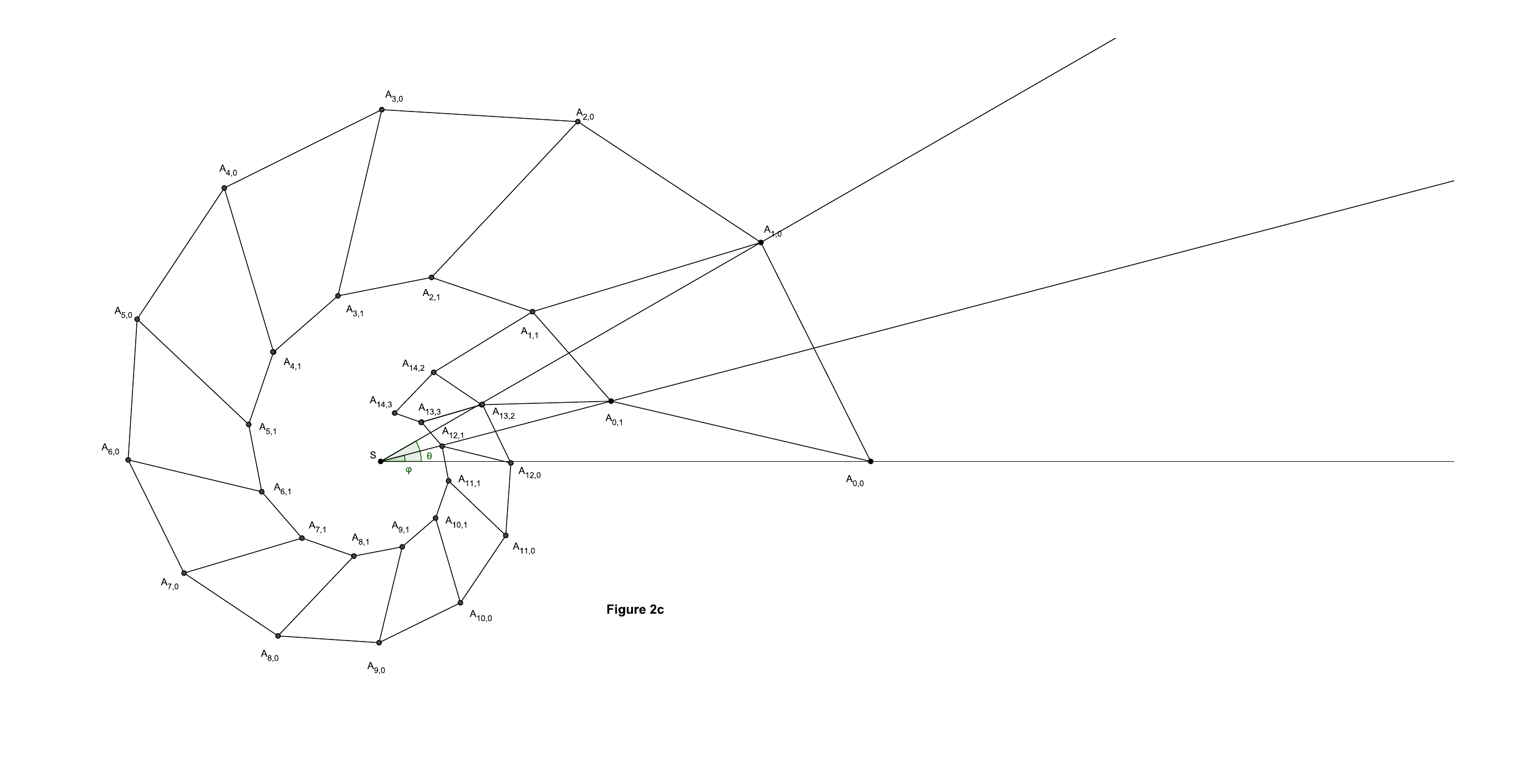}{}
\end{figure}

This is better understood by considering a "travelling" segment which starts from the initial segment's $A_{13,1}A_{0,0}$ position, goes to $A_{14,1}$ (first step) in the $\lambda$ direction, then to $A_{15,1}$ (second step) and so on, covering in 13 steps all positions in order to "arrive" finally to $A_{26,2}$ (which forms an angle $\widehat{\phi}$ with the segment $A_{13,1}A_{0,0}$). This "travelling" segment rotates at a total angle $13\widehat{\theta}$ in order to arrive at the final position. In the same way another "travelling" segment starts from the same initial segment $A_{13,1}A_{0,0}$ in the $\kappa$ direction and goes to $A_{26,2}A_{13,1}$ in one step, rotating at an angle of $\widehat{\phi}$. These two "traveling" segments in order to meet together and appear as one segment, apart from rotating in order to have the same final orientation, they have to have the same size. Because of the branch $B\kappa_{0}$ and the ratio of similarity ratio $\kappa$, we have:\\

$A_{26,2}A_{13,1}=\kappa A_{13,1}$,\\

Also, due to the fact that any quadrangle found in the $\lambda$ direction has a ratio of similarity $\lambda$ with the two adjacent quadrangles found in the same direction and because the first "travelling" segment "travels" in 13 steps in that direction, we have:\\

$A_{26,2}A_{13,1}=\lambda^{13} A_{13,1}$\\

\begin{figure}[bht]
\includegraphics{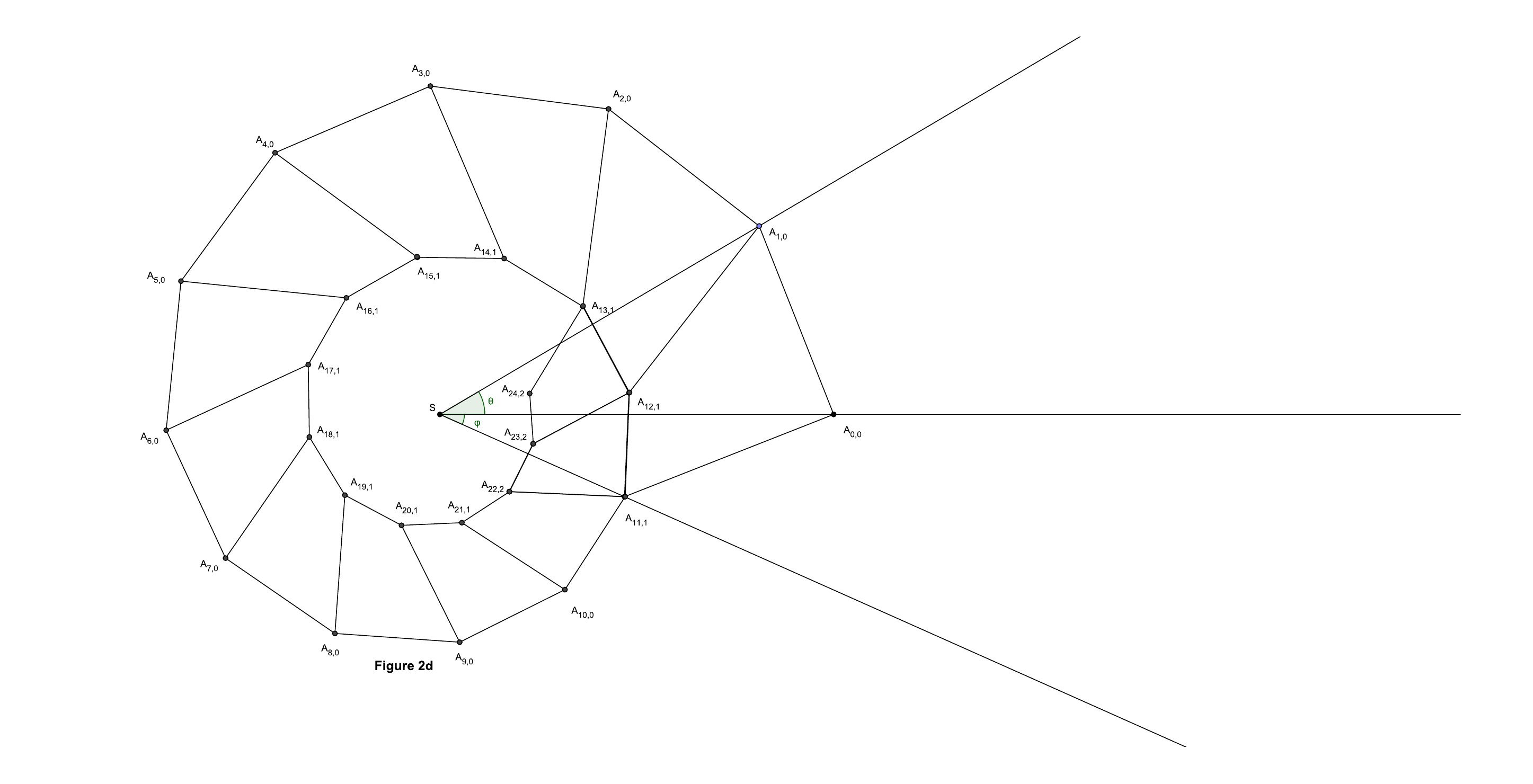}{}
\end{figure}

Apart from the issues of size and rotation, there is also another fact which makes the two "travelling" segments to meet together. These two segments at their final position have the same orientation, the same size and additionally have to form with the center $S$ of the spiral a triangle with a known size and orientation of one side ($A_{26,2}A_{13,1}$, their size when they meet) and angles $\widehat{\phi}$, $\widehat{\sigma}$, $\pi-\widehat{\phi}-\widehat{\sigma}$ so there is only one triangle satisfying these requirements, so they meet together and appear as one segment.

The two equations above give us

\begin{equation}
 \lambda=\sqrt[13]{\kappa}
\end{equation}

The equations deduced from the example of figure 2b let us do the following: given the defining parameters which are $\widehat{\phi}$, $\kappa$ plus an initial segment such as $A_{0,0}A_{13,1}$ and the conditions that there is only one branch $B\lambda_{0}$ and thirteen branches $B\kappa{i}$ ($i=0,..,12$), we can construct a closed quadrangle spiral system.\\

\textbf{RULE 3}. In a closed spiral system where the number of branches is finite, we define as $n$, the parameter for the number of branches in the $\lambda$ direction and $m$, the parameter for the number of branches in the $\kappa$ direction.\\

In the case of figure 2b, we have $n=13$ and $m=1$. In figure 2c we have the same defining parameters as in figure 2b, except for $m=2$, which is the number of branches in the $\kappa$ direction. Thus, the initial branch $B\lambda_{0}$ meets the initial branch $B\kappa_{0}$ at vertex $A_{13,2}$ ($A_{i,j}$, where $i=0+13,j=0+2$) after their common start and in a similar way the second branch $B\lambda_{1}$ and the initial branch $B\kappa_{0}$ meet at vertex $A_{13,3}$ ($A_{i,j}$, where $i=0+13,j=1+2$) after their common start at vertex $A_{0,1}$. The branches $B\lambda_{0}$ and $B\kappa_{0}$ meet again at point $A_{26,4}$ ($A_{i,j}$, where $i=0+2\times13,j=0+2\times2$), the branches $B\lambda_{1}$ and $B\kappa_{0}$ at point $A_{26,5}$ ($A_{i,j}$, where $i=0+2\times13,j=1+2\times2$), and so on. This is how the vertices numbering system of the closed spiral systems works.\\

The equations 3 and 4 become respectively: $13\widehat{\theta}=2\pi+2\widehat{\phi}$ and $\lambda=\sqrt[13]{\kappa^{2}}$. Their general form is:

\begin{equation}
n\widehat{\theta}-m\widehat{\phi}=2\pi  \Leftrightarrow n\widehat{\theta}=2\pi+m\widehat{\phi}
\end{equation}

\begin{equation}
\lambda^{n}=\kappa^{m} \Leftrightarrow \lambda=\sqrt[n]{\kappa^{m}}
\end{equation}

It has to be noted that the values of $\kappa$ and $\lambda$ which are related to equation 6 have to be either greater or less than one (usually $\kappa<1$, so we have to have $\lambda<1$ as well), so the equations 1, 2, 8, 17 have to be treated accordingly (if the branches $B\kappa_{i}$ and $B\lambda_{i}$ are contra-rotating, the reciprocal expressions of either $\kappa$ or $\lambda$ have to be used).\\

In figure 2d we have a spiral system of which the branches of $\kappa$ and $\lambda$ directions contra-rotate, where $n=11$ and $m=1$. The only difference is in the equation 5 (general form) which becomes:

\begin{equation}
n\widehat{\theta}+m\widehat{\phi}=2\pi \Leftrightarrow n\widehat{\theta}=2\pi-m\widehat{\phi}
\end{equation}

The above equation states that the "traveling" segment in figure 2d (which starts from the initial segment's $A_{13,1}A_{0,0}$ position, as in figure 2b) rotates covering a total angle $11\widehat{\theta}$ ($n=11$) which is equal to $2\pi$ minus the angle $\widehat{\phi}$, (since $\widehat{\phi}$ is the angle between the segments $A_{0,0}A_{11,1}$ and $A_{11,1}A_{22,2}$ and $m=1$). The equation 7 is necessary but not sufficient in order to have a closed quadrangular spiral system with contra-rotating branches, so the right combination of conditions of RULE 2 has to apply, which in this case is:\\

$\widehat{\sigma}<\pi/2-\widehat{\phi}/2$ and $(\widehat{\omega}+\widehat{\sigma})>\pi/2-\widehat{\theta}/2$ or\\

$\widehat{\sigma}>\pi/2-\widehat{\phi}/2$ and $(\widehat{\omega}+\widehat{\sigma})<\pi/2-\widehat{\theta}/2$\\

The following parameters belong to an example where equation 7 holds but is not sufficient, $\widehat{\phi}=30 degr$, $n=6$, $m=4$ and $\kappa=0.75$, from equation 7 we get: $\widehat{\theta}=40 degr$, from equation 1 we get: $\widehat{\sigma}=46.94 degr$, from equation 6 we get: $\lambda=0.8256$  and from equation 2: $\widehat{\omega}=8.353 degr$. From these data and the conditions of RULE 2 stated above, it is obvious that the spiral system does not have contra-rotating branches and therefore it cannot be closed.\\

In order to "design" a closed spiral system, we specify the basic parameters $\widehat{\phi}$, $\kappa$ (or $\widehat{\theta}$, $\lambda$), $n$ and $m$ and we calculate $\widehat{\theta}$ (or $\widehat{\phi}$) from equation 5 or 7 (co-rotating or contra-rotating), $\lambda$ (or $\kappa$) from equation 6, $\widehat{\sigma}$ from equation 1 and $\widehat{\omega}$ from equation 2. The parameters to be specified can be considered as "design" parameters of the closed spiral system we want to produce. The vertex numbering system in the present work covers all cases of spiral systems, closed or open, with divergence angle (section 6) or not.\\

In Sushida, Hizume, Yamagishi [2, p.2] we notice that the parameters $n$ and $m$ which have a indirect but similar effect on the shape of closed spirals, have to be relatively prime integers (something necessary only when a divergence angle has to be calculated, as shown in section 6), without any geometrical and visual meaning, thus restricting the ability to understand in a realistic way the produced spirals. We also notice that no conditions similar to those of RULE 2 (and the rest of the rules which appear in the present work with their geometrical nature) exist.\\

In [5] a number of examples using the exponential definition of spiral systems is presented and it is quite interesting to see that the governing equations in the complex plane do not have the direct geometrical and visual meaning of such Euclidean methods as the present work.

It is interesting to note that when $m=n$, a condition which creates a symmetrical closed system, we get from equation 5 that:
$\widehat{\phi}=\widehat{\theta}-2\pi/n$. This gives us: $\widehat{\phi}<\widehat{\theta}$  for symmetrical quadrangular or triangular co-rotating closed spiral systems\\

The parameters $n$, $m$, $\widehat{\phi}$, $\widehat{\theta}$, $\kappa$, $\lambda$ of figures 2b, 2c and 2d have values which appear at Table 3 of Appendix 3.

\section{\textbf{CLOSED TRIANGULAR SPIRAL SYSTEMS WHERE $\widehat{\omega}=\widehat{\phi}$}}

As stated previously, any quadrangular spiral system system can become triangular, so a closed triangular system can exist in two cases, first when $\widehat{\omega}=\widehat{\phi}$ (figure 3a) and second when $\widehat{\omega}=0$ (figure 4a). The parameters $n$, $m$, $\widehat{\phi}$, $\widehat{\theta}$, $\kappa$, $\lambda$ of figures 3a, 3e and 4a have values which appear at Table 3 of Appendix 3. In the present section we examine the first case and calculate the relevant necessary conditions, such as the ones which hold for figure 3a, where $n=12$ and $m=2$.

\begin{figure}[bht]
\includegraphics{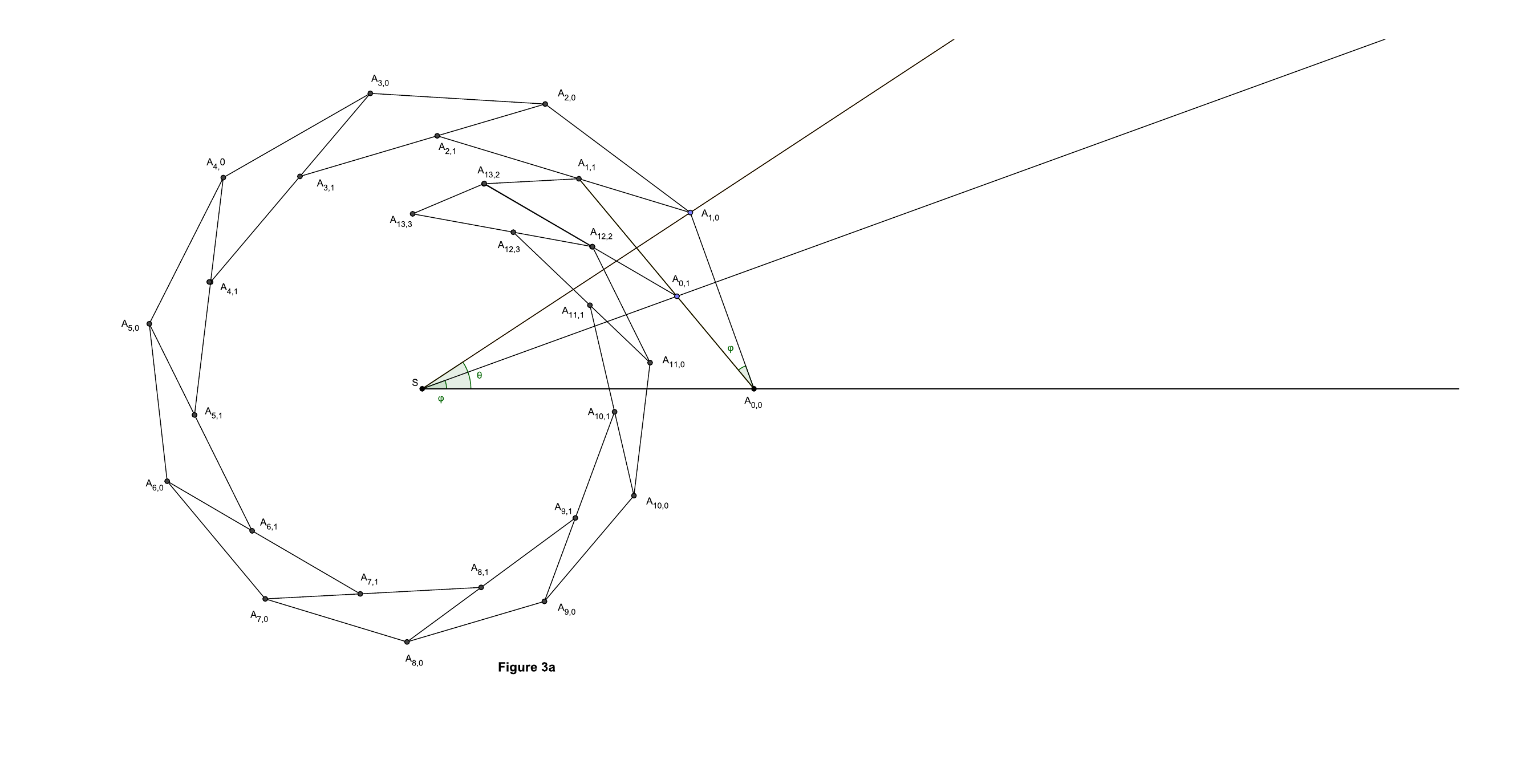}{}
\end{figure}

From equation 2 and since $\widehat{\omega}=\widehat{\phi}$, we have:\\

\begin{eqnarray*}
tan (\widehat{\phi}+\widehat{\sigma})=\lambda sin \widehat{\theta} /(1-\lambda cos \widehat{\theta})\Rightarrow
\end{eqnarray*}

\begin{eqnarray*}
\lambda=(sin\widehat{\phi}+tan\widehat{\sigma} cos\widehat{\phi})/(sin(\widehat{\theta}+\widehat{\phi})+tan\widehat{\sigma}cos(\widehat{\theta}+\widehat{\phi}))
\end{eqnarray*}

\begin{figure}[bht]
\includegraphics{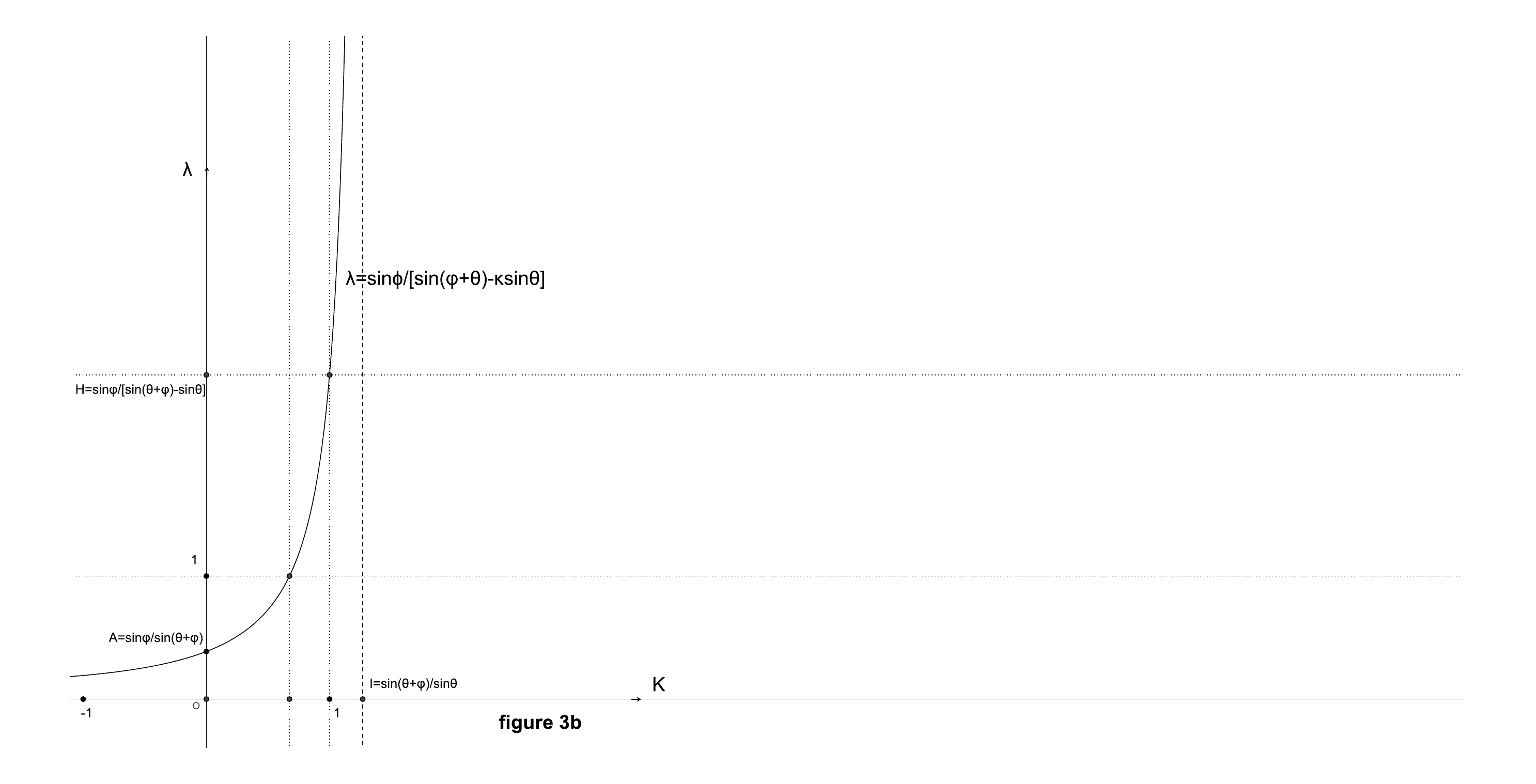}{}
\end{figure}

\begin{figure}[bht]
\includegraphics{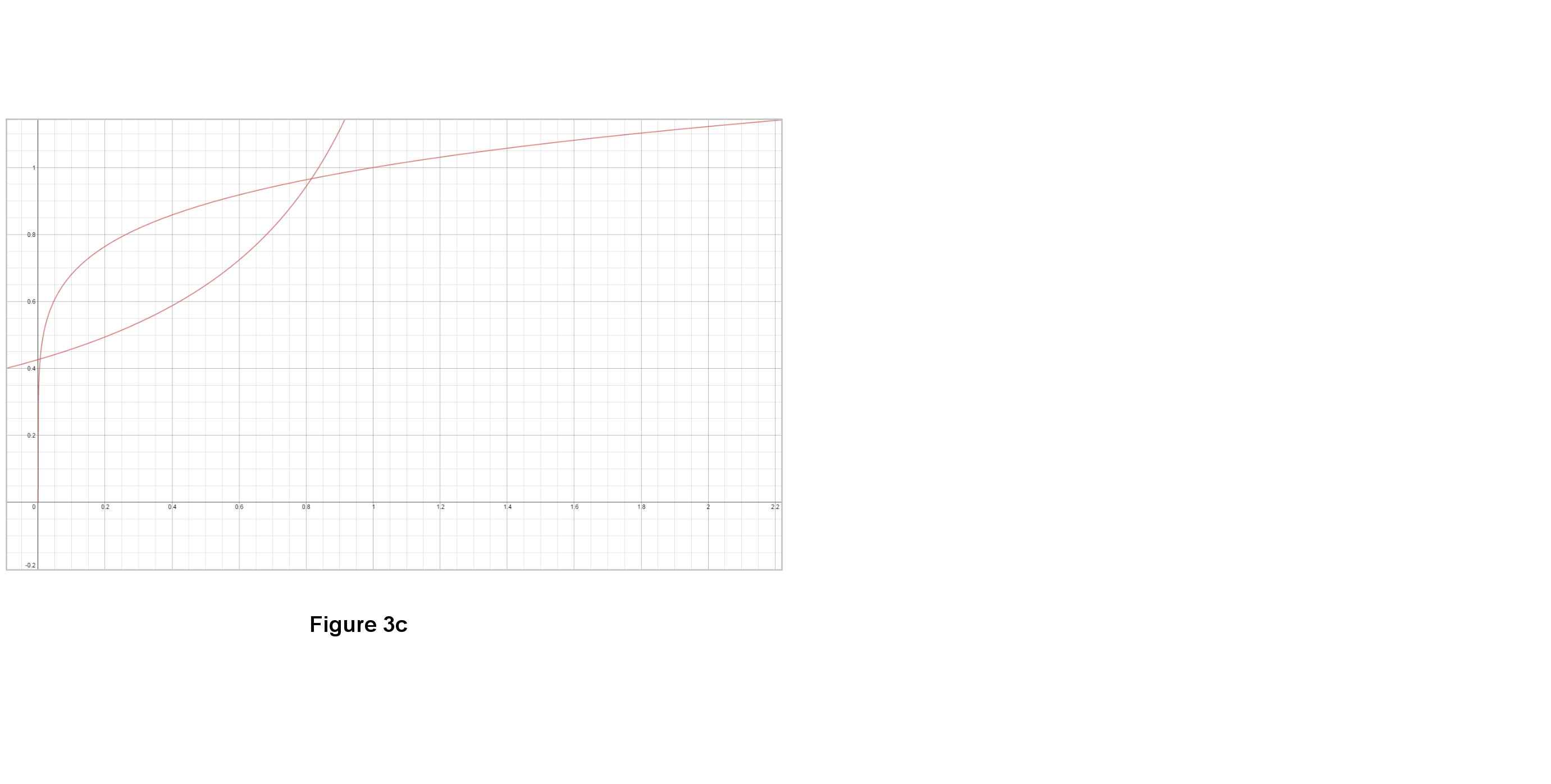}{}
\end{figure}

\begin{figure}[bht]
\includegraphics{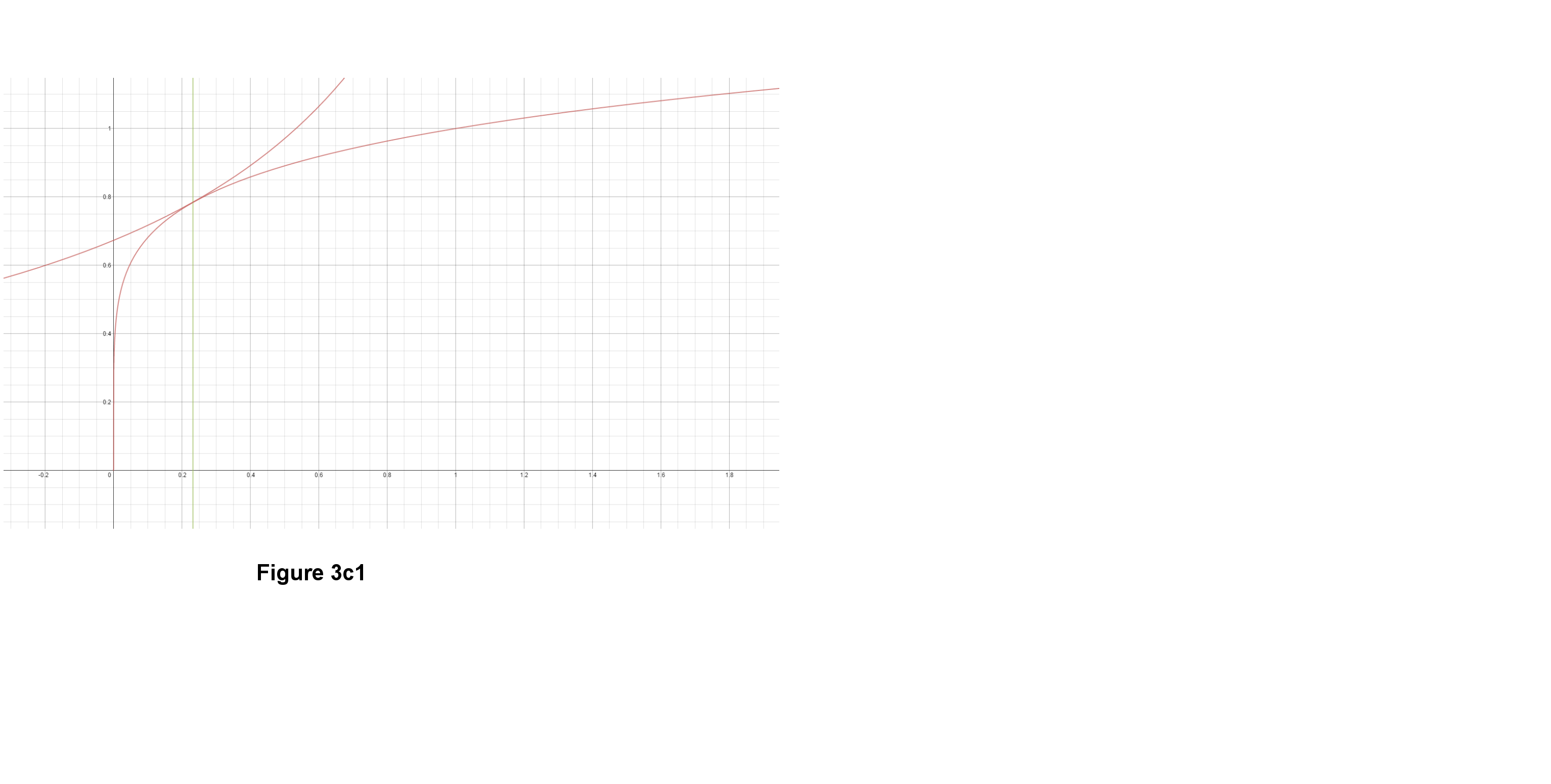}{}
\end{figure}

\begin{figure}[bht]
\includegraphics{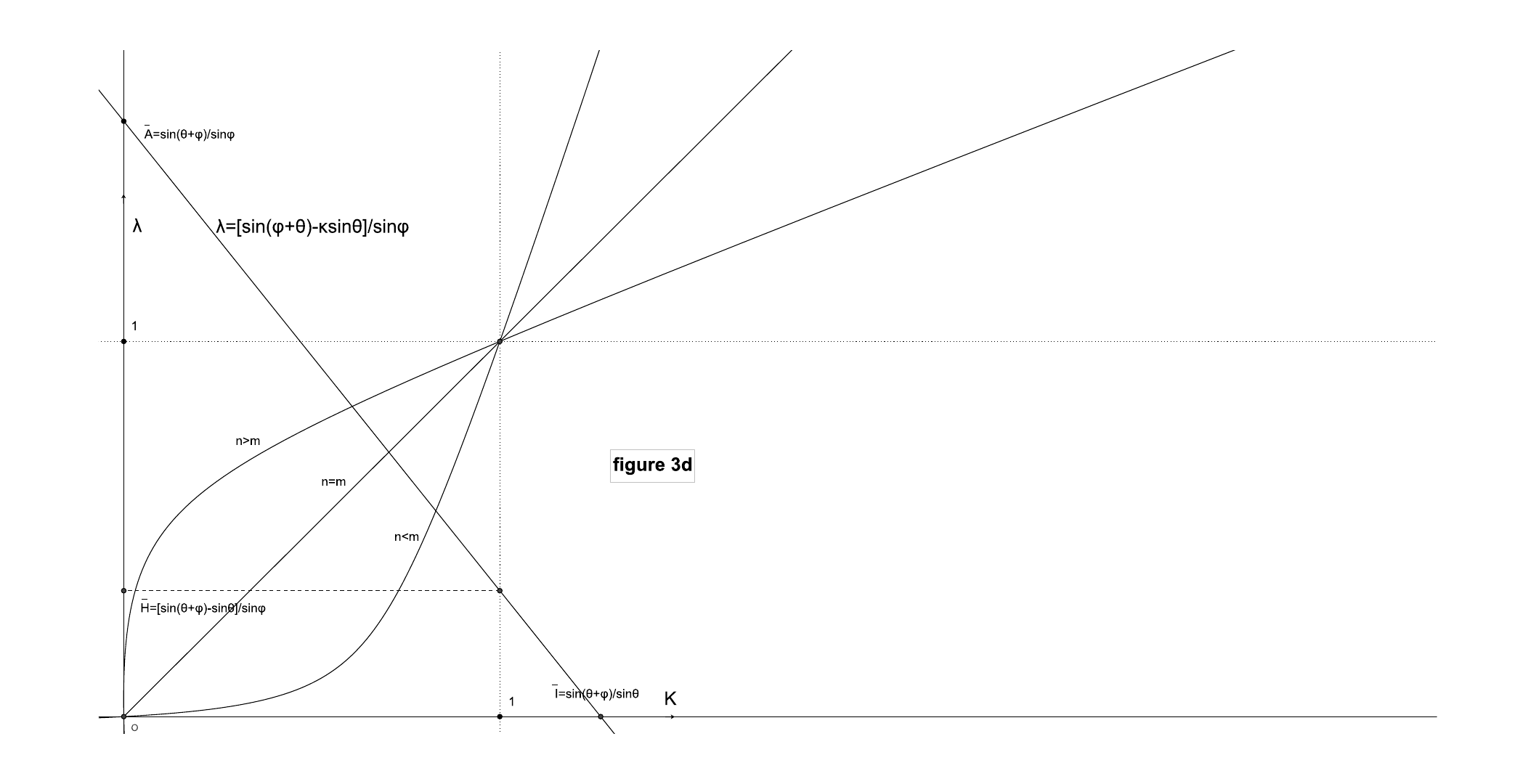}{}
\end{figure}

\begin{figure}[bht]
\includegraphics{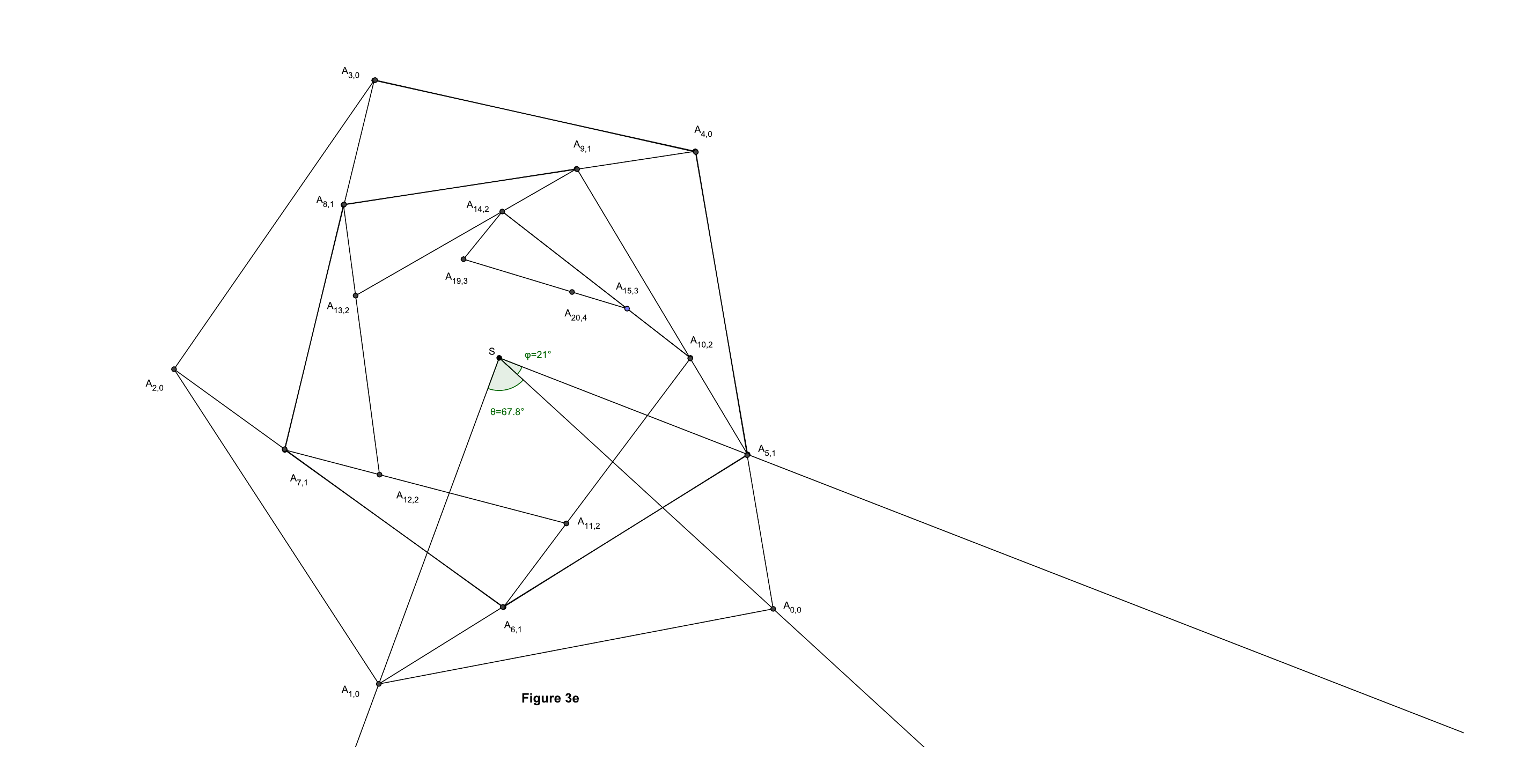}{}
\end{figure}

\begin{figure}[bht]
\includegraphics{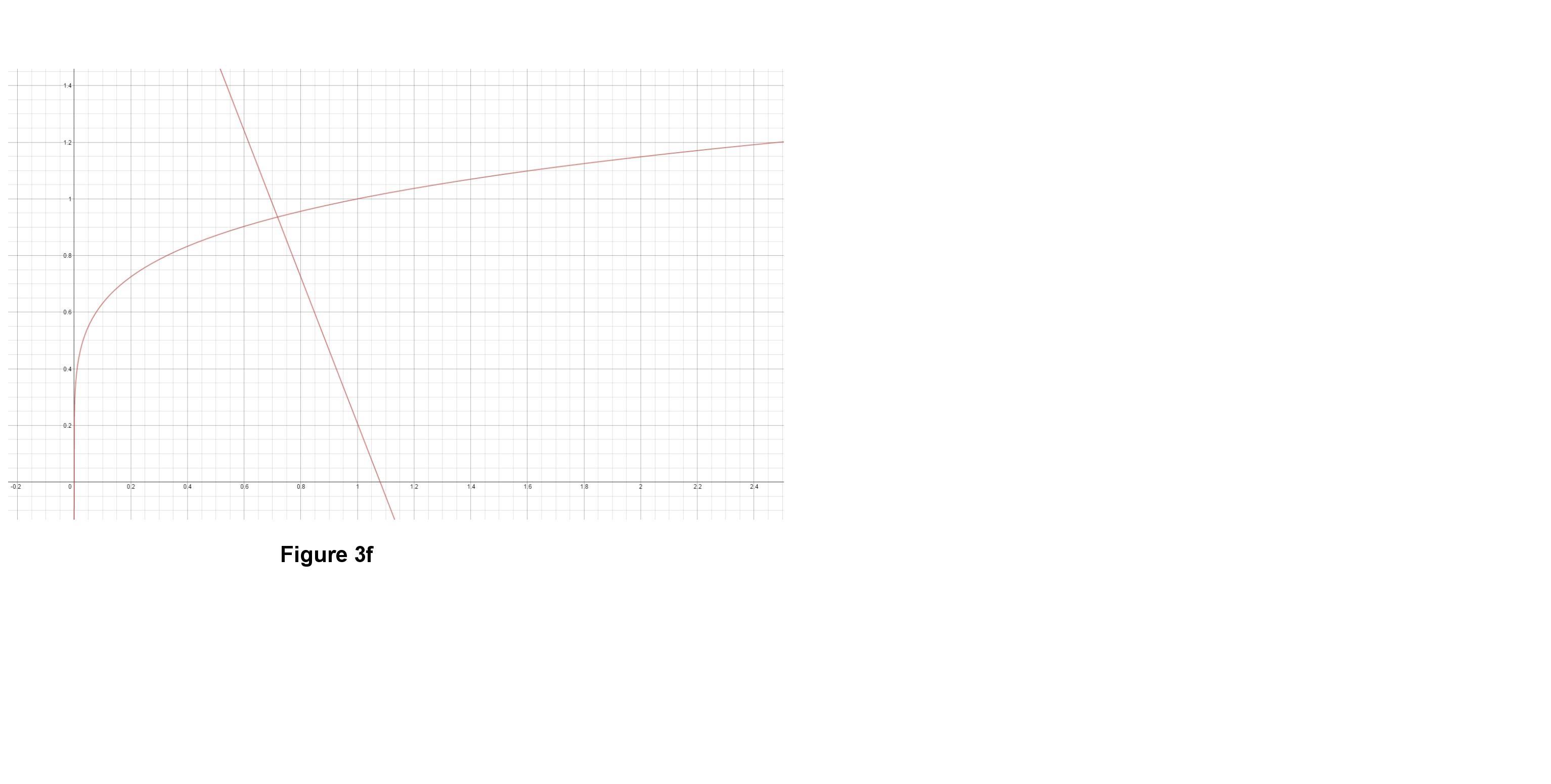}{}
\end{figure}

By using equation 1 in the above equation we get:

\begin{eqnarray*}
\lambda=(sin\widehat{\phi}(1-\kappa cos\widehat{\phi})+\kappa sin\widehat{\phi} cos\widehat{\phi})/ (sin(\widehat{\theta}+\widehat{\phi})(1-\kappa cos\widehat{\phi})+\kappa sin\widehat{\phi}cos(\widehat{\theta}+\widehat{\phi})) \Rightarrow
\end{eqnarray*}

\begin{equation}
\lambda=sin\widehat{\phi}/(sin(\widehat{\theta}+\widehat{\phi})-\kappa sin\widehat{\theta})
\end{equation}

All the parameters of a closed spiral system are $\widehat{\phi}$, $\widehat{\theta}$, $\kappa$,  $\lambda$, $n$, $m$, $\widehat{\sigma}$ and $\widehat{\omega}$. In order to "design" a triangular spiral system, we specify $\widehat{\phi}$, $n$, $m$, we get $\widehat{\theta}$ from equation 5, or 7 (co-rotating or contra-rotating branches) and find a pair of values $\kappa$,  $\lambda$ which satisfy equations 8 and 6, using tools such as GraphSketch or Desmos ($\widehat{\omega}=\widehat{\phi}$ and $\widehat{\sigma}$ is calculated by equation 2). Given the values of $\widehat{\phi}$, $n$, $m$, it is not always possible to find solutions for the values of $\kappa$,  $\lambda$. Considering that the parameters $\kappa$ and $\lambda$ have to have values less than 1 in order to match the equations 6 and 8, we have to identify two cases:\\

\textbf{CASE 1}. Branches $B\kappa_{i}$ and $B\lambda_{i}$  co-rotate\\

The parameters $n$, $m$, $\widehat{\phi}$, $\widehat{\theta}$ are related with equation 5. Given that the branches co-rotate anticlockwise (figure 3a), equation 8 gives values of $\lambda$ less than 1 as it is ($SA_{0,0}>SA_{1,0}$). In figure 3b we have a graph of equation 8 with all the important points of $\kappa$ axis and $\lambda$ axis. Obviously the equation 6 has its graph with concave up when $n<m$ and concave down when $n>m$. By studying the range of values of these points of equation 8 in relation to equation 6, we can identify the range of values of $\widehat{\phi}$, $n$, $m$ which allow us to calculate solutions for the values of $\kappa$,  $\lambda$. Table 1 below gives these points and their definitions.

\begin{center}
\textbf{TABLE 1}\
\end{center}
\renewcommand{\arraystretch}{1.35}
\begin{center}
\begin{tabular}{|c|}
  \hline

  $A=sin\widehat{\phi}/sin(\widehat{\theta}+\widehat{\phi})$\\

  $I=sin(\widehat{\theta}+\widehat{\phi})/sin\widehat{\theta}$\\

  $H=sin\widehat{\phi}/(sin(\widehat{\theta}+\widehat{\phi})-sin\widehat{\theta})$\\

  \hline

\end{tabular}
\end{center}

$ $\\

In figure 3c we have the graphs of equations 6 and 8 for the closed co-rotating triangular spiral system of figure 3a, where we notice that there are two solutions of $\kappa$ and $\lambda$ values satisfying the equations, one with too small $\kappa$ value for any visual presentation and the other, related with the system of figure 3a.

From figure 3b it is obvious that if $A\geq 1$ then there is no solution. This and some other conditions create the following rules which apply to the case of closed triangular spiral system where $\widehat{\omega}=\widehat{\phi}$.\\

\textbf{RULE 4}. If $A\geq 1$ then there are no values of $\kappa$,  $\lambda$ which satisfy the equations 6 and 8. Therefore solutions can be found only when $A<1$, which means $sin\widehat{\phi}<sin(\widehat{\theta}+\widehat{\phi})$ or $\widehat{\theta}+\widehat{\phi}<\pi-\widehat{\phi}$ which is equivalent to:

\begin{equation}
\widehat{\theta}+2\widehat{\phi}<\pi
\end{equation}

Two immediate corollaries of this are that there can be no solution when $\widehat{\phi}=\widehat{\theta}=\pi/3$ (when the spiral triangles become equilateral) and also that:

\begin{equation}
\widehat{\phi}<\pi/2
\end{equation}

\textbf{RULE 5}. If $I\geq 1$ then $sin\widehat{\theta}<sin(\widehat{\theta}+\widehat{\phi})$ or $\widehat{\theta}+\widehat{\phi}<\pi-\widehat{\theta}$ which is equivalent to:

\begin{equation}
\widehat{2\theta}+\widehat{\phi}<\pi
\end{equation}

This inequality together with inequality 9 give us:

\begin{equation}
\widehat{\theta}+\widehat{\phi}<2\pi/3
\end{equation}

\textbf{RULE 6}. If $I<1$ then $sin\widehat{\theta}>sin(\widehat{\theta}+\widehat{\phi})$ or $\widehat{\theta}+\widehat{\phi}>\pi-\widehat{\theta}$ which is equivalent to:

\begin{equation}
\widehat{2\theta}+\widehat{\phi}>\pi
\end{equation}

\textbf{RULE 6a}. Inequalities 9 and 13 give us the following:

\begin{equation}
\widehat{\theta}-\widehat{\phi}>0
\end{equation}

\begin{equation}
(\pi-\widehat{\phi})/2<\widehat{\theta}<\pi-2\widehat{\phi}
\end{equation}

\textbf{RULE 6b}. Since $\widehat{\phi} max=\widehat{\theta}$, from inequality 15 we get:

\begin{equation}
0<\widehat{\phi}<\pi/3
\end{equation}
\\

\textbf{RULE 7}. As proved in Appendix 2, part 1, we always have $H>1$ and when $I\geq 1$, these conditions obviously give us two solutions, provided that Rule 4 holds and the appropriate values of $\widehat{\phi}$, $n$, and $m$ are chosen, such as those in figure 3c, mentioned above. When $I<1$, it is obvious that again we have two solutions, provided that the above mentioned conditions apply. In Appendix 2, part 2 these conditions are examined and also in Table 4, maximum values of $\widehat{\phi}$ (which give us one double solution as defined in Appendix 2 part 2) and their related $\kappa$ values for specific values of $n$ and $m$ are presented.\\

\textbf{CASE 2}. Branches $B\kappa_{i}$ and $B\lambda_{i}$  contra-rotate\\

The parameters $n$, $m$, $\widehat{\phi}$, $\widehat{\theta}$ are related with equation 7. Given that the branches contra-rotate, equation 8 gives values of $\lambda$ greater than 1 (because $SA_{0,0}<SA_{1,0}$, whereas $SA_{0,0}>SA_{0,1}$ and these two conditions guarantee the contra-rotating characteristics, being equivalent to RULE 2 when $\widehat{\omega}=\widehat{\phi}$), so in order to match the equation 6, it has to be raised to the power of -1 and become reciprocal as follows:\\

\begin{eqnarray*}
\lambda=(sin(\widehat{\theta}+\widehat{\phi})-\kappa sin\widehat{\theta})/sin\widehat{\phi}
\end{eqnarray*}

The equivalent variables $A,H,I$ of Table 1, become reciprocal in this case as follows:

\newpage

\begin{center}
\textbf{TABLE 1A}\
\end{center}
\renewcommand{\arraystretch}{1.35}
\begin{center}
\begin{tabular}{|c|}
  \hline

  $\overline{A}=sin(\widehat{\theta}+\widehat{\phi})/sin\widehat{\phi}$\\

  $\overline{I}=sin\widehat{\theta}/sin(\widehat{\theta}+\widehat{\phi})$\\

  $\overline{H}=(sin(\widehat{\theta}+\widehat{\phi})-sin\widehat{\theta})/sin\widehat{\phi}$\\

  \hline

\end{tabular}
\end{center}

$ $\\

In figure 3d we have a graph of the reciprocal of equation 8 which is a straight line, with all the important points of Table 1A, plus the three versions of equation 6, concave up when $n<m$, concave down when $n>m$ and straight line when $n=m$.\\

\textbf{RULE 8}. Since we have $H>1$ as proven in Appendix 2, part 1, it is obvious that $\overline{H}<1$ in this case. This condition guarantees always the existence of one solution for all cases of equation 6, $n<m$, $n=m$, $n>m$, as shown in figure 3d. An example of such a spiral system with contra-rotating branches exists in figure 3e (the values of its parameters appear in Table 3) and the graphs of its relevant equations (the equivalent of those in graph 3d) exist in figure 3f.\\

Contrary to the case of closed quadrangular systems with contra-rotating branches, the equation 7 is sufficient together with the reciprocal version of equation 8 which takes into account the necessary condition in order to have values of $\lambda$ less than 1, as mentioned previously at the beginning of CASE 2 paragraph. The only limits for the angles are :$\widehat{\phi}<2\pi/m$ and $\widehat{\theta}<2\pi/n$ from equation 7.

\section{\textbf{CLOSED TRIANGULAR SPIRAL SYSTEMS WHERE $\widehat{\omega}=0$}}

As in the previous section, we examine and calculate the necessary conditions to have a closed triangular spiral system in this case, such as the one of figure 4a. The parameters $n$, $m$, $\widehat{\phi}$, $\widehat{\theta}$, $\kappa$, $\lambda$ of figure 4a have values which appear at Table 3 of Appendix 3. From equations 1 and 2 and for $\widehat{\omega}=0$ we have:\\

\begin{figure}[bht]
\includegraphics{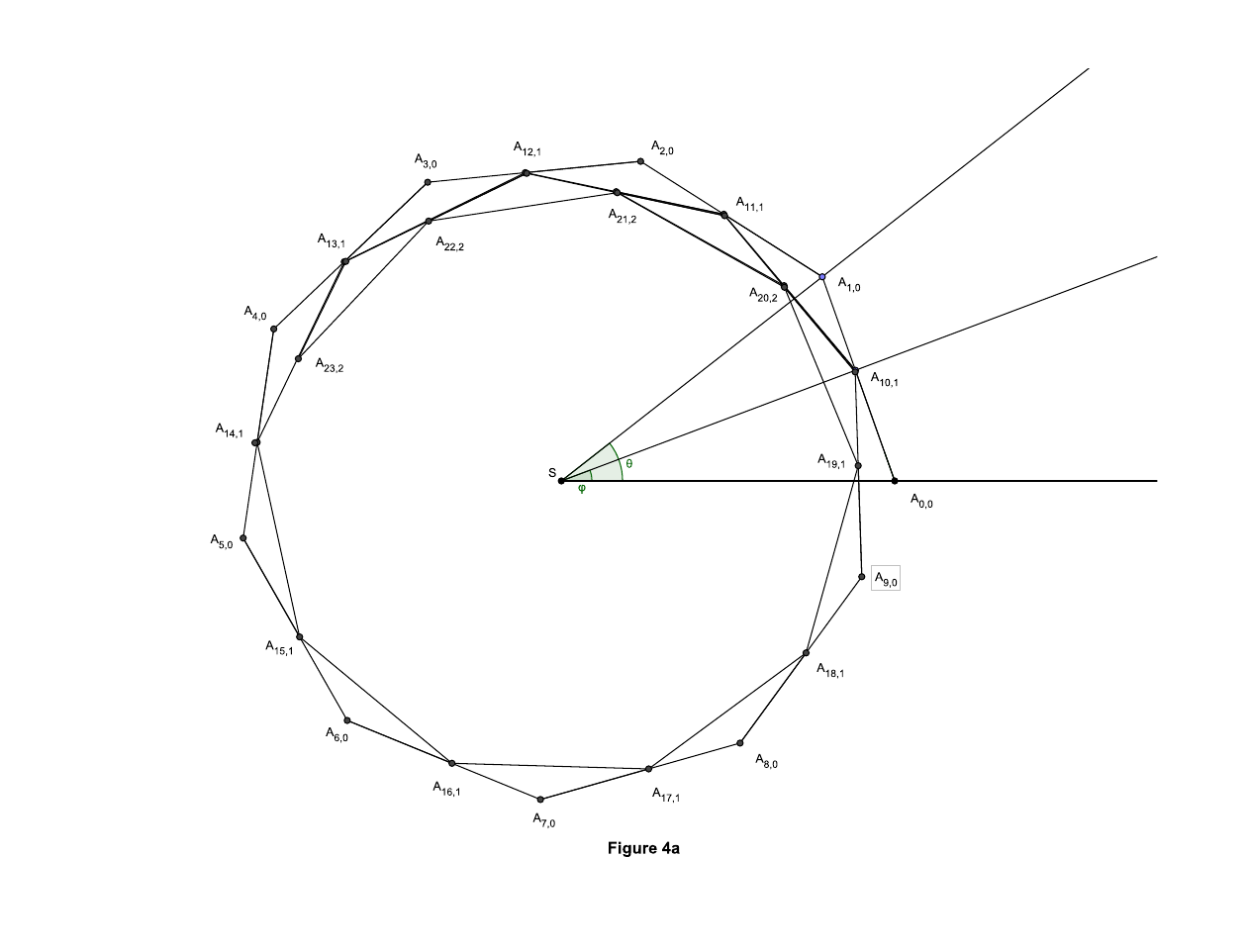}{}
\end{figure}

\begin{figure}[bht]
\includegraphics{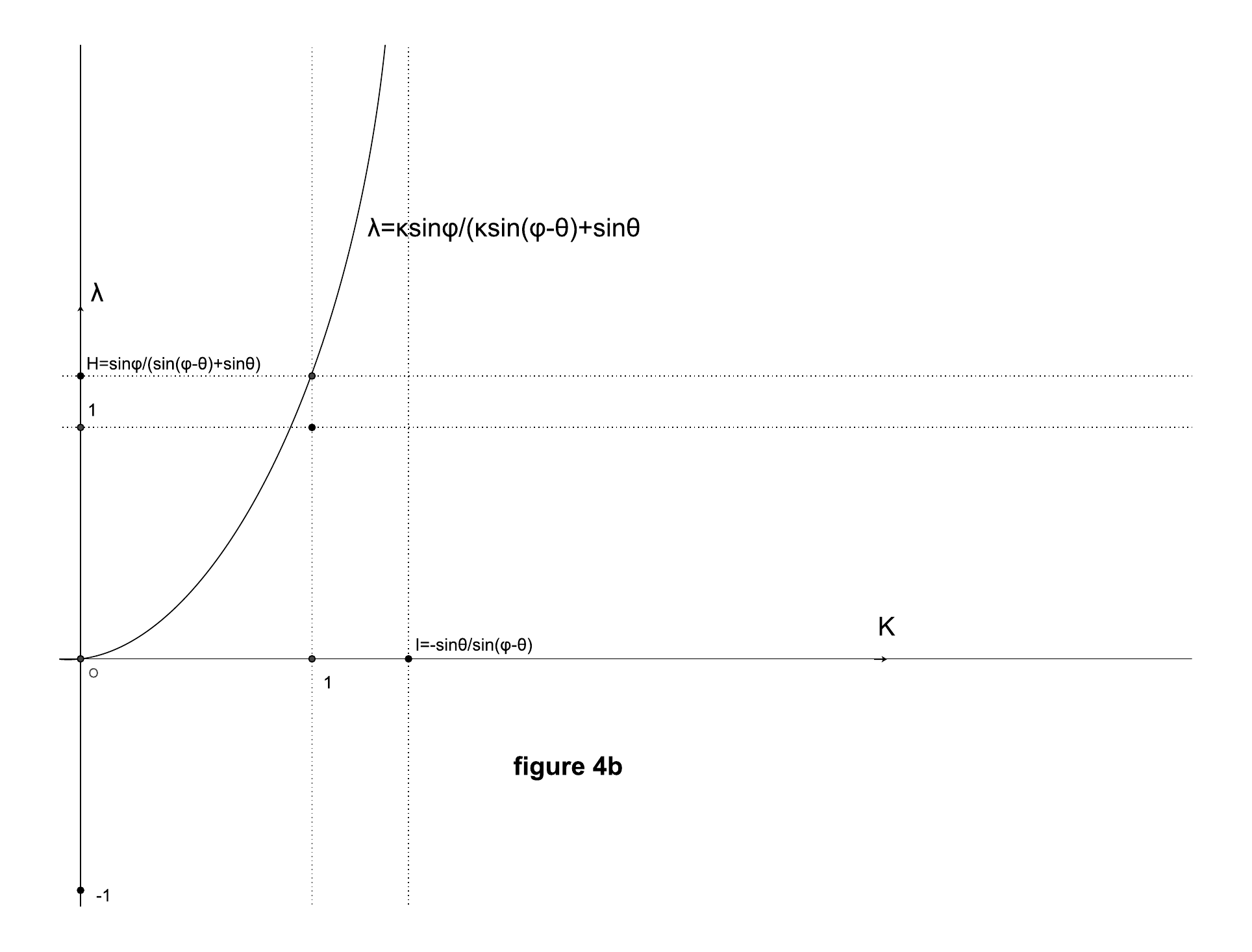}{}
\end{figure}

\begin{figure}[bht]
\includegraphics{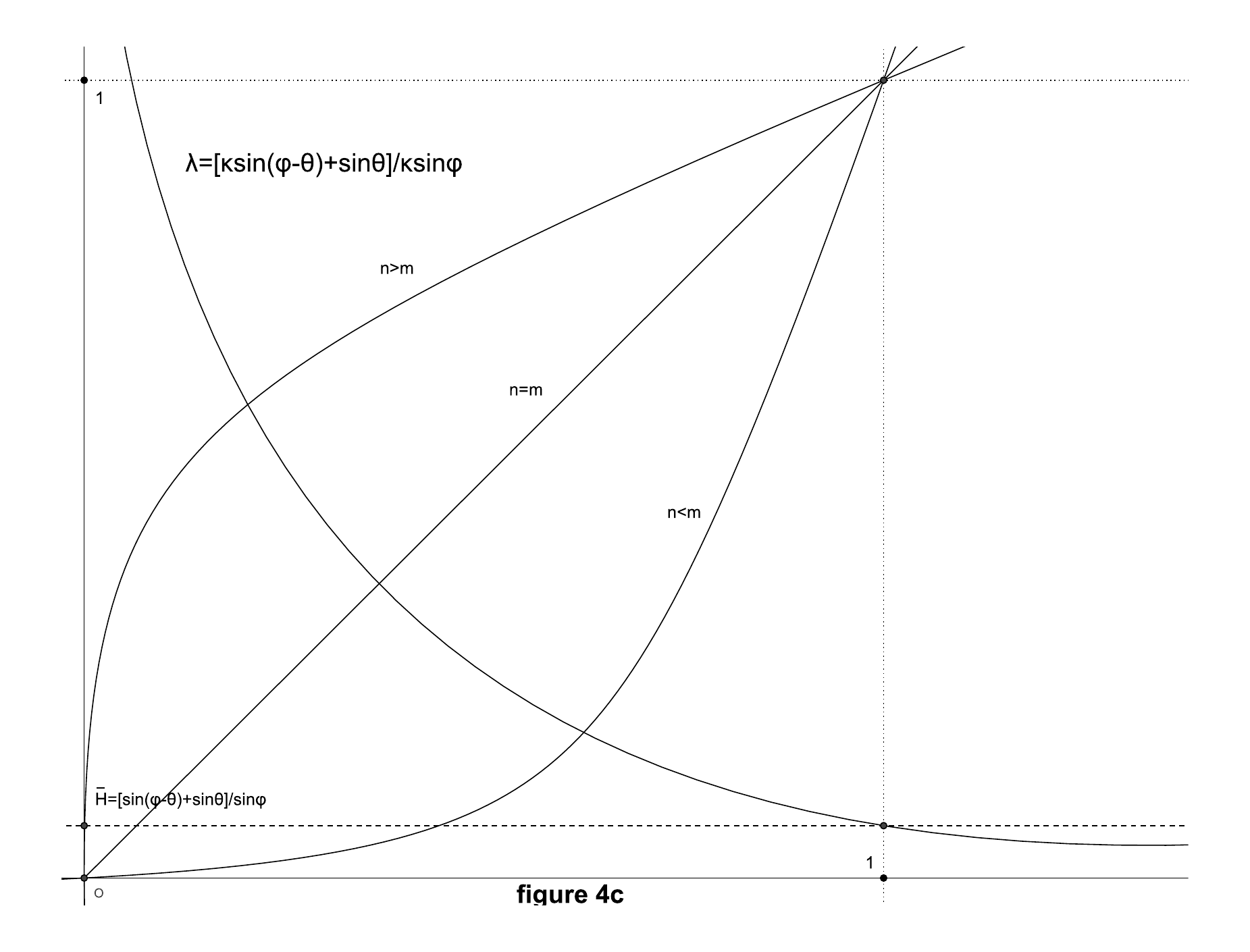}{}
\end{figure}

\begin{figure}[bht]
\includegraphics{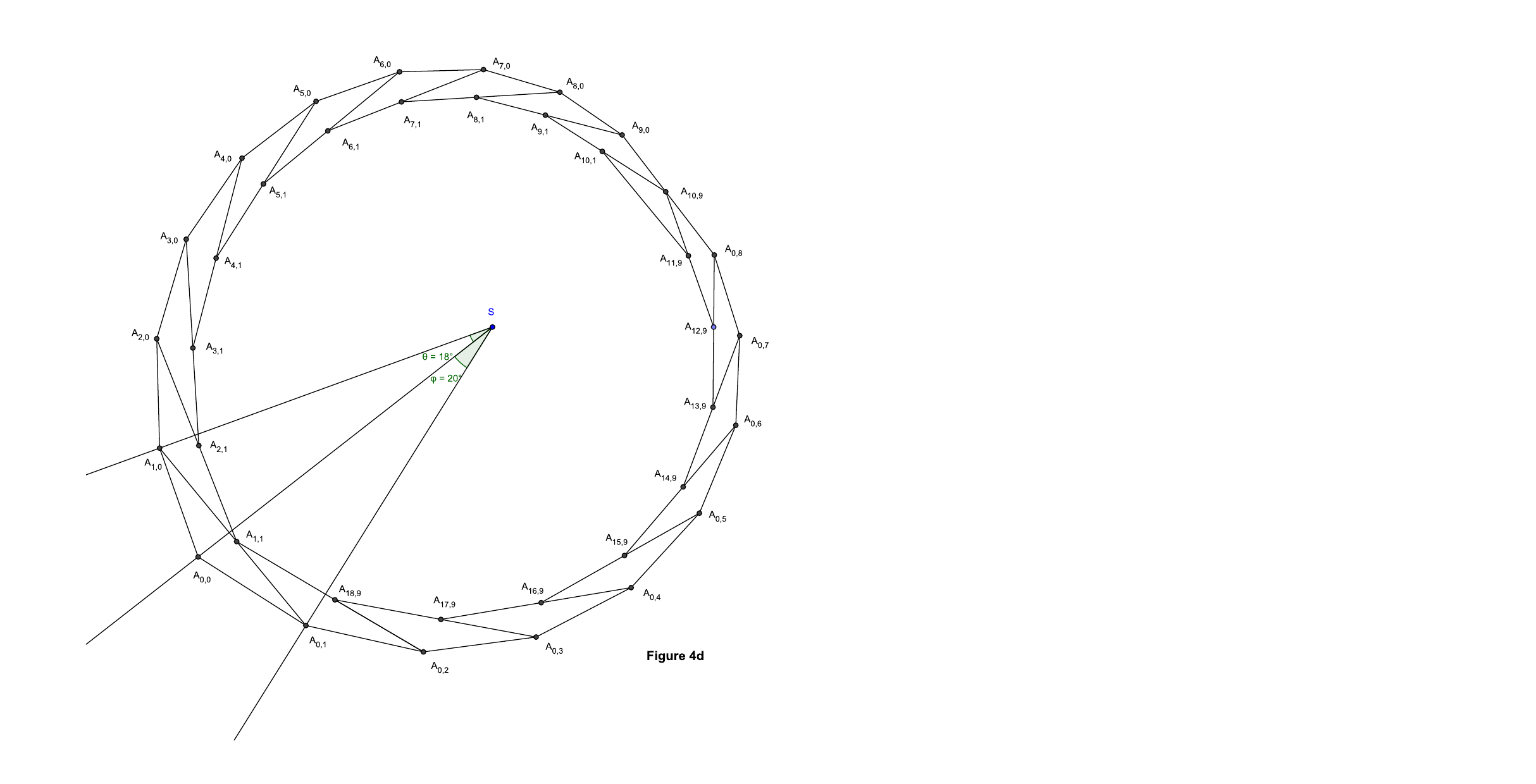}{}
\end{figure}

$\kappa sin\widehat{\phi}/(1-\kappa cos\widehat{\phi})=\lambda sin\widehat{\theta}/(1-\lambda cos\widehat{\theta})\Rightarrow$

\begin{equation}
 \lambda=\kappa sin\widehat{\phi}/(sin\widehat{\theta}+\kappa sin(\widehat{\phi}-\widehat{\theta}))
\end{equation}

As in the previous section, we have to distinguish two cases :\\

\textbf{CASE 1}. Branches $B\kappa_{i}$ and $B\lambda_{i}$  co-rotate\\

The parameters $n$, $m$, $\widehat{\phi}$, $\widehat{\theta}$ are related with equation 5. In figure 4b, we have a graph of equation 17 with all the important points of $\kappa$ and $\lambda$ axis. Table 2 below gives these points and their definitions.

\begin{center}
\textbf{TABLE 2}\
\end{center}
\renewcommand{\arraystretch}{1.35}
\begin{center}
\begin{tabular}{|c|}
  \hline

  $I=-sin\widehat{\theta}/sin(\widehat{\phi}-\widehat{\theta})$\\

  $H=sin\widehat{\phi}/(sin(\widehat{\phi}-\widehat{\theta})+sin\widehat{\theta})$\\

  \hline
\end{tabular}
\end{center}

$ $\\

As in the previous section, by studying the range of values of these points of equation 17 in relation to equation 6, we can identify the conditions related to values of $\widehat{\phi}$, $n$, $m$ which allow us to calculate solutions for the values of $\kappa$,  $\lambda$. These conditions create the following rules which apply to the case of closed triangular spiral system where $\widehat{\omega}=0$.\\

\textbf{RULE 9}. As proved in Appendix 2, part 2, when $0<\widehat{\phi}<\widehat{\theta}$ we have $H>1$ and these conditions always give us two solutions as long as $n>m$, one of which (for $\kappa=0$) is without geometrical significance, as shown shown in figure 4c in which the graphs of equations 6 and 17 are presented for the case of figure 4a.\\

\textbf{RULE 9a}. When $0<\widehat{\phi}<\widehat{\theta}$ and if $I>1$, then $sin\widehat{\theta}>sin(\widehat{\theta}-\widehat{\phi})$ or $\widehat{\theta}-\widehat{\phi}<\pi-\widehat{\theta}$, which is equivalent to:

\begin{equation}
2\widehat{\theta}-\widehat{\phi}<\pi
\end{equation}

\textbf{RULE 9b}. When $0<\widehat{\phi}<\widehat{\theta}$ and if $I<1$, then :

\begin{equation}
2\widehat{\theta}-\widehat{\phi}>\pi
\end{equation}

\textbf{RULE 10}. When $0<\widehat{\phi}=\widehat{\theta}$ then from equation 17 we get $\lambda=\kappa$, which does not produce any meaningful solution.\\

\textbf{RULE 11}. When $\widehat{\phi}>\widehat{\theta}$ it is obvious from the definition of $I$ that always $I<0$.\\

\textbf{RULE 11a}. When $\widehat{\phi}>\widehat{\theta}$ and $n<m$ it is obvious that we have two solutions, one of which (for $\kappa=0$) is without geometrical significance.\\

\textbf{CASE 2}. Branches $B\kappa_{i}$ and $B\lambda_{i}$  contra-rotate\\

The parameters $n$, $m$, $\widehat{\phi}$, $\widehat{\theta}$ are related with equation 7. Given that the branches contra-rotate, equation 17 gives values of $\lambda$ greater than 1 (because $SA_{0,0}<SA_{1,0}$ whereas $SA_{0,0}>SA_{0,1}$ as in CASE 2 of section 3), so in order to match the equation 6, it has to be raised to the power of -1 and become reciprocal as follows:

\begin{eqnarray*}
 \lambda=(sin\widehat{\theta}+\kappa sin(\widehat{\phi}-\widehat{\theta}))/\kappa sin\widehat{\phi}
\end{eqnarray*}

The equivalent variable $H$ of Table 2, becomes reciprocal in this case as follows:\\

$\overline{H}=(sin(\widehat{\phi}-\widehat{\theta})+sin\widehat{\theta})/sin\widehat{\phi}$\\

In figure 4c we have a graph of the reciprocal of equation 17 with the important point $\overline{H}$, plus the three versions of equation 6, concave up when $n<m$, concave down when $n>m$ and straight line when $n=m$.\\

\textbf{RULE 12}. As proved in Appendix 2, part 2, when $0<\widehat{\phi}<\widehat{\theta}$ we have $H>1$, so $\overline{H}<1$, which gives us always one solution for all cases of equation 6, as shown in figure 4c. When $\widehat{\phi}>\widehat{\theta}$ we have $\overline{H}>1$, which gives us no solution.\\

As in section 3 and contrary to the case of closed quadrangular systems with contra-rotating branches, the equation 7 is sufficient together with the reciprocal version of equation 17 which takes in account the necessary condition in order to have values of $\lambda$ less than 1, as mentioned previously at the beginning of CASE 2 paragraph.\\

\section{\textbf{EQUIVALENCE OF CLOSED TRIANGULAR SPIRAL SYSTEMS WHERE EITHER $\widehat{\omega}=0$ OR $\widehat{\overline{\omega}}=\widehat{\phi}$}}

As stated in the introduction section, we can treat the case of a spiral system where $\widehat{\omega}=0$ just like the case of a new spiral system where $\widehat{\overline{\omega}}=\widehat{\phi}$ by introducing the parameters $\overline{\lambda}$, $\overline{m}$, and the angle $\widehat{\overline{\theta}}$ of the new system, which has the other parameters ($\kappa$, $\widehat{\phi}$ and $n$) the same with the first one. In this section we will examine the necessary relations in order to have the equivalence of the two types of closed triangular systems, described in sections 3 and 4. We have all together five cases of equivalence, presented below. It has to be noted that the analysis and the equations 20 to 32 can be applicable also to any quadrangular and triangular spirals.\\

\textbf{CASE 1}. Branches $B\kappa_{i}$ and $B\lambda_{i}$  co-rotate, $\widehat{\theta}>\widehat{\phi}$ and $n>m$\\

In this case we have equation 5 and equation 6 for the old system as shown in figure 4a, where we have $\overline{\lambda}=A_{20,2}A_{11,1}/A_{11,1}A_{2,0}=\kappa/\lambda$ and  $\widehat{\overline{\theta}}=\widehat{A_{20,2}SA_{11,1}}=\widehat{A_{11,1}SA_{2,0}}=\widehat{\theta}-\widehat{\phi}$ similarly to figure 1d, taking in consideration that the branches $B\kappa_{i}$ and $B\lambda_{i}$ co-rotate, $m<n$ and $\kappa<\lambda$.\\

In Appendix 4 we calculate the equation 20 and it is obvious that the relevant equations of the other cases can be calculated in a similar way. In figure 4d we have the equivalent new system  ($\widehat{\overline{\omega}}=\widehat{\phi}$) to the old system ($\widehat{\omega}=0$) related to figure 4a. The parameters $n$, $\overline{m}$, $\widehat{\phi}$, $\widehat{\overline{\theta}}$, $\kappa$, $\overline{\lambda}$ of figure 4d have values which appear at Table 3 of Appendix 3. The following equations apply for both systems:

\begin{equation}
 n\widehat{\overline{\theta}}=2\pi-\overline{m}\widehat{\phi}=2\pi-(n-m)\widehat{\phi}
\end{equation}

\begin{equation}
\overline{\lambda}=\kappa/\lambda
\end{equation}

\begin{equation}
\overline{\lambda}=\sqrt[n]{\kappa^{\overline{m}}}=\sqrt[n]{\kappa^{n-m}}=\kappa^{(n-m)/n}
\end{equation}

At the equation 20 we have that $\widehat{\overline{\theta}}=\widehat{\theta}-\widehat{\phi}$ and $\overline{m}=n-m$, therefore equations 5 and 20 are equivalent, also \textbf{the branches of the new system contra-rotate}. By supposing that $\kappa<1$ and because $n>m$ and equations 6 and 21 hold, we have that $\overline{\lambda}<1$ and $\kappa<\lambda$. Finally from equations 6, 21 and 22 we have:\\

$\overline{\lambda}=\kappa/\lambda=\kappa^{(n-m)/n}$ and $\lambda=\sqrt[n]{\kappa^{m}}=\kappa^{m/n}\Rightarrow\kappa=\kappa^{m/n}\kappa^{(n-m)/n}=\kappa$\\

From the above we deduce that all the equations of the old and the new system are equivalent (being of the same type) and compatible.\\

\textbf{CASE 2}. Branches $B\kappa_{i}$ and $B\lambda_{i}$  co-rotate, $\widehat{\theta}>\widehat{\phi}$ and $n<m$\\

In this case we have equation 5 and equation 6 for the old system, plus the following equations for the new system:

\begin{equation}
 n\widehat{\overline{\theta}}=2\pi+\overline{m}\widehat{\phi}=2\pi+(m-n)\widehat{\phi}
\end{equation}

\begin{equation}
\overline{\lambda}=\lambda/\kappa
\end{equation}

\begin{equation}
\overline{\lambda}=\sqrt[n]{\kappa^{\overline{m}}}=\sqrt[n]{\kappa^{m-n}}=\kappa^{(m-n)/n}
\end{equation}

The branches of the new system co-rotate. By supposing that $\kappa<1$ and because $n<m$ and equations 6 and 24 hold, we have that $\overline{\lambda}<1$ and $\kappa>\lambda$. From equations 6, 24 and 25 we have:\\

$\overline{\lambda}=\lambda/\kappa=\kappa^{(m-n)/n}$ and $\lambda=\sqrt[n]{\kappa^{m}}=\kappa^{m/n}\Rightarrow1/\kappa=\kappa^{-m/n}\kappa^{(m-n)/n}=1/\kappa$\\

\textbf{CASE 3}. Branches $B\kappa_{i}$ and $B\lambda_{i}$  co-rotate, $\widehat{\theta}<\widehat{\phi}$ and $n<m$\\

In this case we have the following equation (a version of equation 5) and equation 6 for the old system :

\begin{equation}
m\widehat{\phi}=2\pi+n\widehat{\theta}
\end{equation}

We also have the following equations for the new system:

\begin{equation}
 n\widehat{\overline{\theta}}=2\pi-\overline{m}\widehat{\phi}=2\pi-(m-n)\widehat{\phi}
\end{equation}

\begin{equation}
\overline{\lambda}=\lambda/\kappa
\end{equation}

\begin{equation}
\overline{\lambda}=\sqrt[n]{\kappa^{\overline{m}}}=\sqrt[n]{\kappa^{m-n}}=\kappa^{(m-n)/n}
\end{equation}

The branches of the new system contra-rotate. By supposing that $\kappa<1$ and because $n<m$ and equations 6 and 28 hold, we have that $\overline{\lambda}<1$ and $\kappa>\lambda$. Finally from equations 6, 28 and 29 we have:\\

$\overline{\lambda}=\lambda/\kappa=\kappa^{(m-n)/n}$ and $\lambda=\sqrt[n]{\kappa^{m}}=\kappa^{m/n}\Rightarrow1/\kappa=\kappa^{-m/n}\kappa^{(m-n)/n}=1/\kappa$\\

\textbf{CASE 4}. Branches $B\kappa_{i}$ and $B\lambda_{i}$  co-rotate, $\widehat{\theta}<\widehat{\phi}$ and $n>m$\\

In this case we have equation 26 and equation 6 for the old system, plus the following equations for the new system:

\begin{equation}
 n\widehat{\overline{\theta}}=2\pi+\overline{m}\widehat{\phi}=2\pi+(n-m)\widehat{\phi}
\end{equation}

\begin{equation}
\overline{\lambda}=\kappa/\lambda
\end{equation}

\begin{equation}
\overline{\lambda}=\sqrt[n]{\kappa^{\overline{m}}}=\sqrt[n]{\kappa^{n-m}}=\kappa^{(n-m)/n}
\end{equation}

The branches of the new system co-rotate. By supposing that $\kappa<1$ and because $n>m$ and equations 6 and 31 hold, we have that $\overline{\lambda}<1$ and $\kappa<\lambda$. From equations 6, 31 and 32 we have:\\

$\overline{\lambda}=\kappa/\lambda=\kappa^{(n-m)/n}$ and $\lambda=\sqrt[n]{\kappa^{m}}=\kappa^{m/n}\Rightarrow\kappa=\kappa^{m/n}\kappa^{(n-m)/n}=\kappa$\\

\textbf{CASE 5}. Branches $B\kappa_{i}$ and $B\lambda_{i}$  contra-rotate\\

In this case we have equation 7 and equation 6 for the old system, plus the following equations (equation 34 is calculated at Appendix 1.5) for the new system:

\begin{equation}
 n\widehat{\overline{\theta}}=2\pi-\overline{m}\widehat{\phi}=2\pi-(n+m)\widehat{\phi}
\end{equation}

\begin{equation}
\overline{\lambda}=\kappa\lambda
\end{equation}

\begin{equation}
\overline{\lambda}=\sqrt[n]{\kappa^{\overline{m}}}=\sqrt[n]{\kappa^{n+m}}=\kappa^{(n+m)/n}
\end{equation}

The branches of the new system contra-rotate. By supposing that $\kappa<1$ and $\lambda<1$, we have that $\overline{\lambda}<1$. From equations 6, 34 and 35 we have:\\

$\overline{\lambda}=\kappa\lambda=\kappa^{(n+m)/n}$ and $\lambda=\sqrt[n]{\kappa^{m}}=\kappa^{m/n}\Rightarrow\kappa=\kappa^{-m/n}\kappa^{(n+m)/n}=\kappa$\\

The parameters listed below belong to an example of an old system and a new equivalent spiral system related to this case, following the conditions of RULE 2: $m=2$, $n=8$, $\widehat{\phi}=15 degr$, $\widehat{\theta}=41.25 degr$, $\kappa=0.945$, $\lambda=0.986$, $\widehat{\sigma}=70.38 degr$, $\overline{\lambda}=0.932$, $\widehat{\overline{\theta}}=26.25 degr$, $\overline{m}=10$.

\section{\textbf{Divergence Angle}}

The divergence angle $\widehat{d}$, as described in [2, p.2 and p.3] (figure 1, angle formed by vertex 1, centre of spiral O and vertex 2), in [3, p.100] in [4, p.20] and in figure 5 (present work, see further below), is defined by drawing lines from two successive primordia to the center of the spiral and $\widehat{d}>\pi$ or $\widehat{d}<\pi$ according to the direction of rotation the angle is measured. As it is stated in [2, p.19] (Bravais-Bravais theorem), given any starting vertex (primordium) at a spiral system with parastichies $(m,n)$, the next vertex found in the $m$ parastichy passing from the starting vertex (in the $\lambda$ direction, according to the present work) differs by $m$ (from the starting vertex) and similarly the next vertex found in the $n$ parastichy passing from the starting vertex (in the $\kappa$ direction) differs by $n$ (from the starting vertex, example in [2, p2]). This means that starting from any vertex of the spiral tiling-system, one has to turn $m\widehat{d}$ times to go to the next vertex of the $m$ parastichy passing from the starting point and  $n\widehat{d}$ times to go to the next vertex of the $n$ parastichy passing from the starting point. This rule can be expressed by the following obvious equations, applicable in the spiral systems described in the present work, when the branches of $\kappa$ and $\lambda$ directions ($m$ and $n$ parastichies according to [4]) co-rotate:

\begin{figure}[bht]
\includegraphics{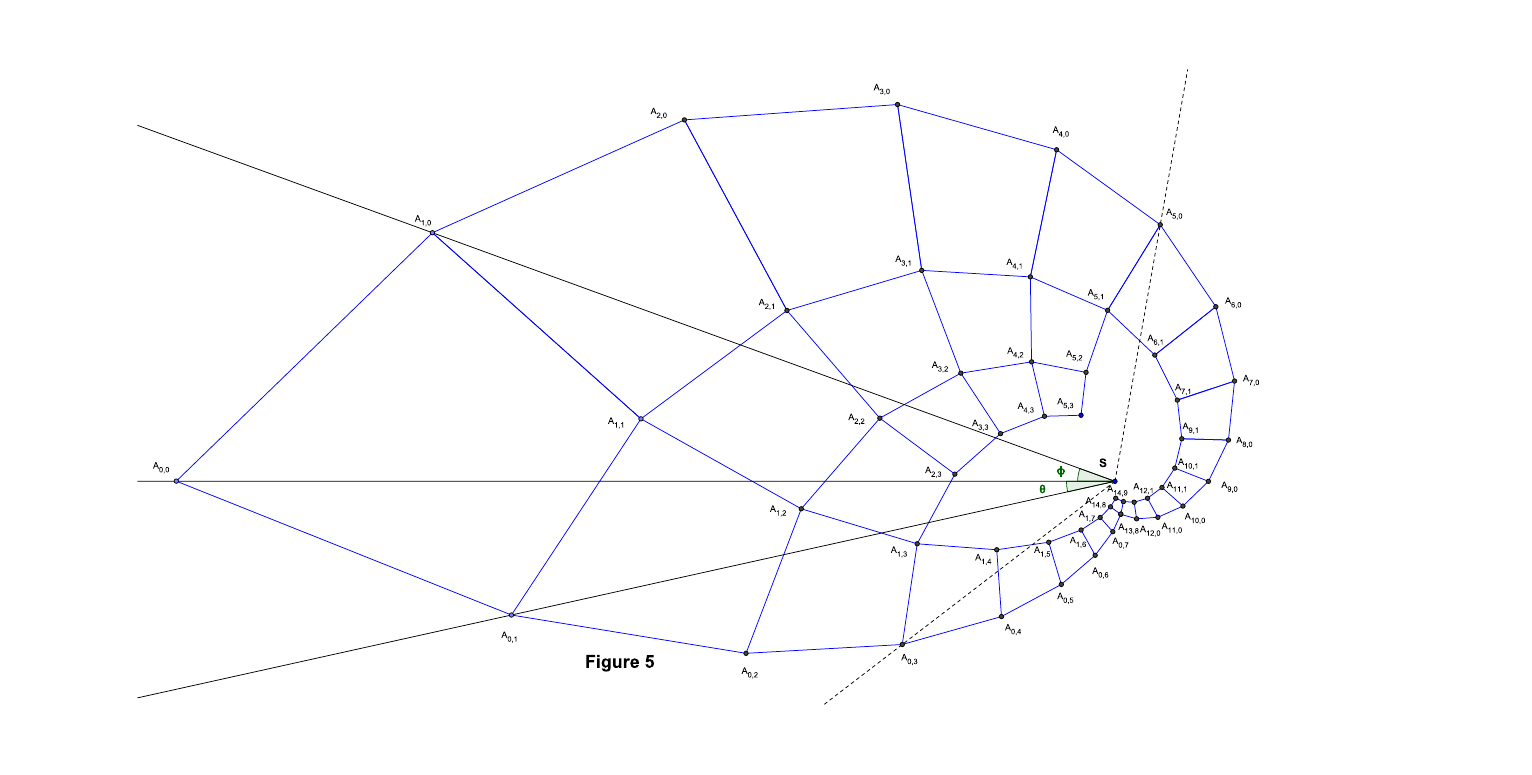}{}
\end{figure}

\begin{equation}
2\pi I_{\lambda}-\widehat{\phi}=n\widehat{d}
\end{equation}

\begin{equation}
2\pi I_{\kappa}-\widehat{\theta}=m\widehat{d}
\end{equation}

Similarly, when the branches of $\kappa$ and $\lambda$ directions contra-rotate, equation 37 remains and equation 36 becomes:

\begin{equation}
2\pi I_{\lambda}+\widehat{\phi}=n\widehat{d}
\end{equation}

In the above equations, $I_{\lambda}$ and $I_{\kappa}$ represent the integer number of $2\pi$ or 360 degr. turns which are required for the rotation to go to the next vertex of each branch, as described above (for examples, see figure 5 and further below). If the opposite rotation is to be applied, the above equations change sign and different values $I'_\lambda$, $I'_\kappa$ and $\widehat{d'}$ are used, as shown in Appendix 5, part 2.\\
From the above equations, the following ones (equation 39 when branches co-rotate, equation 40 when branches contra-rotate) can be deduced, as shown in Appendix 5 part 1:

\begin{equation}
\widehat{d}=I_{\lambda} \widehat{\theta}-I_{\kappa} \widehat{\phi}
\end{equation}

\begin{equation}
\widehat{d}=I_{\lambda} \widehat{\theta}+I_{\kappa} \widehat{\phi}
\end{equation}

In Table 3 the values of the parameters $I_\kappa$, $I_\lambda$ and $\widehat{d}$ have been calculated for the figures of the present work. It is worth mentioning that at figure 3a where $n=12$ and $m=2$, therefore non relatively prime, we have no calculation of $\widehat{d}$. Also in the example of figure 2b ($m=1$, $n=13$, $\widehat{d}=331.2degr.$), if the point $A_{13,1}$ is chosen as starting point, then going clockwise $I_\lambda=12$ times with angle increments of $\widehat{\theta}=28.8degr.$ we arrive at point $A_{1,0}$ and from there we go anticlockwise with $I_\kappa=1$ increment of $\widehat{\phi}=14.6degr.$, we arrive at point $A_{14,1}$, implementing equations 37 and 39 once ($m=1$), in order to get what we stated at the beginning of this section, "one has to turn $m\widehat{d}$ times to go to the next vertex of the $m$ parastichy passing from the starting point". Applying this process another 12 times, therefore a total of 13 times ($n=13$), we arrive at point $A_{26,2}$, implementing equation 36, in order to have turned "$n\widehat{d}$ times to go to the next vertex of the $n$ parastichy passing from the starting point". The case of figure 5 is based on [2,p2, contra-rotating example of figure 1]. The relevant parameters are, $m=13$, $n=8$, $\widehat{\theta}=12.5degr.$, $\widehat{\phi}=20degr.$, $\widehat{d}=137.5degr.$, $\kappa=0.78$, $\lambda=0.668$, $I_\lambda=3$, $I_\kappa=5$ (the last two parameters were calculated by using equation 41), which all comply with the above mentioned equations. In figure 5, assuming that the starting vertex is $A_{0,3}$, we have that $\widehat{d}=\widehat{A_{0,3}SA_{5,0}}=137.5degr.$ It is interesting to notice that the ray starting from point $S$ and passing from point $A_{0,3}$ (first primordium) goes through points $A_{0,3}\rightarrow A_{1,3}\rightarrow A_{2,3}\rightarrow A_{3,3}\rightarrow A_{4,3}\rightarrow A_{5,3}$, thus turning $5\times 20 degr.$ ($\kappa$ direction) and continues going through points $A_{5,3}\rightarrow A_{5,2}\rightarrow A_{5,1}\rightarrow A_{5,0}$, thus turning another $3\times12.5degr.$ ($\lambda$ direction), in order to pass from point $A_{5,0}$ (second primordium). This is a visual application of equation 40. The visual applications of equations 38 and 37 are produced by assuming that this ray (starting from $A_{0,3}S$) has to rotate $n\widehat{d}$ $(8 \times 137.5degr.$) to go $A_{1,3}S$ ($\widehat{A_{0,3}SA_{1,3}}=\phi=20degr.$, equation 38) and $m\widehat{d}$ $(13 \times 137.5degr.$ to go to $A_{0,4}S$, ($\widehat{A_{0,3}SA_{0,4}}=\widehat{\theta}=12.5degr.$, equation 37).\\

The following Diophantine equation is deduced from equations 36, 37, 5 when branches co-rotate or 37, 38 and 7 when branches contra-rotate (Appendix 5, part 2, where it is also shown that the sign + becomes - when alternative versions of equations 36, 37, 38, opposite rotation, $I'_\lambda$, $I'_\kappa$ and $\widehat{d'}$ are used):

\begin{equation}
n I_{\kappa} - m I_{\lambda} =+1
\end{equation}

It is important to note the following: The equations 36, 37 (branches co-rotating), or 37, 38 (branches contra-rotating) have 3 unknowns, $I_{\kappa}$, $I_{\lambda}$ and $\widehat{d}$, plus the condition that $I_{\kappa}$ and $I_{\lambda}$ are integers. Let $n=\nu q$ and $m=\nu p$ (where $\nu$, $p$, $q$, integers plus $p$ and $q$ relatively prime integers), in other words, $n$ and $m$ non relatively prime integers. From the above and equation 41, we get:\\

$qI_{\kappa}-pI_{\lambda}=1/\nu$ This equation is not valid, since $I_{\kappa},I_{\lambda},p,q,\nu$ are integers. Therefore, the following is proven:\\

\textbf{RULE 13}. At any closed spiral system, if $m$ and $n$ are not relatively prime integers, there is no possibility for a convergence angle to exist.\\

Equation 41, can always be solved as a Diophantine equation, where $n$ and $m$ are relatively prime integers.

Finally, according to [4, p20], [2, p3], the plastochrone $R$ ratio is the ratio of distances between two successively numbered primordia and the center of the spiral and $R>1$. This ratio can be defined according to the present work, using equation 39 or 40 in order to go from any primordium to the next one, after a rotation of $\widehat{d}$ (see figure 5). This definition is as follows, as shown in Appendix 5, part 3:

\begin{equation}
r=(1/\lambda)^{I_\lambda}\kappa^{I_\kappa}=\kappa^{1/n}=\lambda^{1/m}<1\rightarrow R=1/r>1
\end{equation}

In figure 5, according to the definition of plastochrone ratio, $R=A_{0,3}S/A_{5,0}S$ and $R=1/r$, therefore we also have a visual application of equation 42 ($r=(1/0.668)^3(0.78)^5=(0.78)^{1/8}=0.97$).\\

In Appendix 5, part 4, an equivalence is shown between specific equations of the present work and of equations at [2].\\

All the equations produced in the present work do not have any approximation nature, since they fully complywith the used analytical and geometrical approach (in [4], some formulae which are quite similar to those in the present work appear to have an approximation nature, for example in [4, p.37 and p.43]).

\section{\textbf{FINAL REMARKS}}

The purpose of this work is to enable the analysis and design of closed spiral systems using a simple Euclidean geometry approach, offering thus a full visual understanding of the parameters which govern these systems. This approach can be enhanced by further research which can produce some practical applications in phyllotaxis problems which could be combined with theorems or methods based in the already existing complex plane theories. Last but not least, one can design as many variants of spiral systems as he or she can imagine by tuning the appropriate design parameters and this can give pleasure to the designer.

\newpage

\section{Appendices}

\textbf{Appendix 1}

\textbf{Part 1.} In figure 1a, let the ray $A_{0,0}A_{0,1}$ pass from a random point $C$ not belonging to the segment $A_{0,0}A_{0,1}$ and similarly the ray $SA_{0,1}$ intersect the segment $A_{0,0}A_{1,0}$ at point $B$. Let $\widehat{\phi'}=\widehat{A_{0,2}A_{0,1}C}$. Since the triangles $\triangle$$A_{0,0}A_{0,1}S$ and $\triangle$$A_{0,1}A_{0,2}S$ are similar and due to the fact that $\widehat{SA_{0,1}A_{0,2}}+\widehat{A_{0,2}A_{0,1}C}+\widehat{CA_{0,1}B}=\pi$, as a result  $\widehat{\phi'}=\widehat{\phi}$ and $\widehat{A_{0,0}A_{0,1}A_{0,2}}=\pi-\widehat{\phi}$, so:\\

$\widehat{(A_{0,0}A_{0,1},A_{0,1}A_{0,2})}=\widehat{\phi}$.\\

\textbf{Part 2.} By applying the sin law in the similar triangles $\triangle$$A_{0,0}A_{0,1}S$ and $\triangle$$A_{0,1}A_{0,2}S$ we get:
$SA_{0,0}/sin\widehat{(A_{0,0}A_{0,1}S)}=SA_{0,1}/sin\widehat{\sigma}$  therefore
$sin\widehat{(A_{0,0}A_{0,1}S)}=sin\widehat{\sigma}/\kappa$ since $\kappa=SA_{0,1}/SA_{0,0}$. Taking in account that $\pi-\widehat{(A_{0,0}A_{0,1}S)}=\widehat{\phi}+\widehat{\sigma}$, so $sin\widehat{(A_{0,0}A_{0,1}S)}=sin \widehat{\phi} cos\widehat{\sigma}+cos \widehat{\phi} sin\widehat{\sigma}$, we deduce from the above:\\

$tan \widehat{\sigma}=\kappa sin \widehat{\phi} /(1-\kappa cos \widehat{\varphi})$.\\

\textbf{Part 3.} From figure 1a, we get $\widehat{(A_{0,0}A_{0,1},A_{1,0}A_{1,1})}=\widehat{A_{0,0}DA_{1,0}}$ also $\widehat{A_{0,0}A_{1,0}D}=\widehat{A_{0,0}A_{1,0}S}+\widehat{SA_{1,0}D}$ and $\widehat{A_{1,0}A_{0,0}D=\widehat{\omega}}$. From $\triangle A_{0,0}A_{1,0}D$ we get $\widehat{A_{0,0}DA_{1,0}}=\pi-\widehat{A_{1,0}A_{0,0}D}-\widehat{A_{00}A_{1,0}D}$. From $\triangle A_{0,0}A_{1,0}S$ we get $\widehat{A_{0,0}A_{1,0}S}=\pi-\widehat{\sigma}-\widehat{\omega}-\widehat{\theta}$ and from segment $A_{1,0}A_{1,1}$ of branch $B\kappa_{1}$ we get $\widehat{SA_{1,0}A_{1,1}}=\widehat{\sigma}$, therefore $\widehat{A_{0,0}A_{1,0}D}=\pi-\widehat{\theta}-\widehat{\omega}$. From all the above, we get:\\

$\widehat{(A_{0,0}A_{0,1},A_{1,0}A_{1,1})}=\widehat{\theta}$\\

\textbf{Part 4.} From figure 1a and because $\widehat{A_{2,0}A_{1,0}A_{1,1}}=\widehat{\omega}$ we get that:\\

$\widehat{A_{0,0}A_{1,0}A_{1,1}}=\pi-\widehat{\theta}-\widehat{\omega}$ and
$\widehat{A_{0,0}A_{1,0}S}=\pi-\widehat{\theta}-\widehat{\omega}-\widehat{\sigma}$\\

In the same figure, since $\widehat{A_{1,1}A_{0,1}A_{0,2}}=\widehat{\omega}$, we get that:\\

$\widehat{A_{0,0}A_{0,1}A_{1,1}}=\pi+\widehat{\phi}-\widehat{\omega}$\\

 In the same figure, we also have that $\widehat{A_{0,1}A_{1,1}A_{1,0}}=\pi-\widehat{\phi}-\widehat{A_{0,1}A_{1,1}A_{1,2}}$ and because $\widehat{A_{0,1}A_{1,1}A_{1,2}}=\pi-\widehat{\theta}-\widehat{\omega}$, this gives us:\\

$\widehat{A_{0,1}A_{1,1}A_{1,0}}=\widehat{\theta}+\widehat{\omega}-\widehat{\phi}$\\

In figure 1c where $\widehat{\omega}=\widehat{\phi}$, the above relations become:\\

$\widehat{A_{0,0}A_{1,0}A_{1,1}}=\pi-\widehat{\theta}-\widehat{\phi}$ plus  $\widehat{A_{0,0}A_{0,1}A_{1,1}}=\pi$ and $\widehat{A_{0,1}A_{1,1}A_{1,0}}=\widehat{\theta}$\\

In figure 1d where $\widehat{\omega}=0$, the above relations become:\\

$\widehat{A_{0,0}A_{1,0}A_{1,1}}=\pi-\widehat{\theta}$ plus $\widehat{A_{0,0}A_{0,1}A_{1,1}}=\pi+\widehat{\phi}$ and
$\widehat{A_{0,1}A_{1,1}A_{1,0}}=\widehat{\theta}-\widehat{\phi}$\\

\textbf{Part 5.} From figure 1d (applicable also to figures 1a, 1b, 1c) we have: $\lambda=A_{1,1}A_{2,0}/A_{0,1}A_{1,0}$ and $\kappa=A_{0,2}A_{1,1}/A_{0,1}A_{1,0}$, as ratios of diagonals to the relevant similar quadrangles (when $\widehat{\omega}\neq0$) which in this case become triangles. Also from this figure we get that the branches $B\kappa_{i}$ and $B\lambda_{i}$ are co-rotating and $\kappa<\lambda$. From these two, we get:\\

$\overline{\lambda}=A_{0,2}A_{1,1}/A_{1,1}A_{2,0}=\kappa/\lambda$\\

Also from the same figure, since $\widehat{A_{0,1}SA_{1,0}}=\widehat{\theta}-\widehat{\phi}$, due to similarity of triangles $\triangle$$A_{0,1}A_{1,0}S$, $\triangle$$A_{0,2}A_{1,1}S$ and $\triangle$$A_{1,1}A_{2,0}S$, we get:\\

$\widehat{\overline{\theta}}=\widehat{A_{0,2}SA_{1,1}}=\widehat{A_{1,1}SA_{2,0}}=\widehat{\theta}-\widehat{\phi}$\\

In the case where the branches $B\kappa_{i}$ and $B\lambda_{i}$ are contra-rotating,  we have that $1/\lambda=A_{1,1}A_{2,0}/A_{0,1}A_{1,0}$ and $\kappa=A_{0,2}A_{1,1}/A_{0,1}A_{1,0}$, so it can be deduced from the above that:\\

$\overline{\lambda}=A_{0,2}A_{1,1}/A_{1,1}A_{2,0}=\kappa\lambda$\\

\textbf{Appendix 2}

\textbf{Part 1.} From Table 1 we have that $H=sin\widehat{\phi}/(sin(\widehat{\theta}+\widehat{\phi})-sin\widehat{\theta})$, which is equation 8 for $k=1$. If we consider $\widehat{\phi}$ as a variable and the rest of the entities as constants with given value, then $H=f(\widehat{\phi})$, where $\widehat{\phi}<\pi$. The first derivative of the above function is as follows:\\

$f'(\widehat{\phi})=sin\widehat{\theta}(1-cos\widehat{\phi})/(sin(\widehat{\theta}+\widehat{\phi})-sin\widehat{\theta})^2>0$ and $0<\widehat{\phi}<\pi$\\

For $\widehat{\phi}=0$ and following Hopital's rule, we have the following:\\

$lim(f(\widehat{\phi}))=[sin\widehat{\phi}]'/[(sin(\widehat{\theta}+\widehat{\phi})-sin\widehat{\theta})]'=1/cos\widehat{\theta}>1$ where $\widehat{\phi}\rightarrow0$,\\

therefore always $H>1$ under the conditions of Rule 7.\\

\textbf{Part 2.} Let $A=f(\widehat{\phi})=sin\widehat{\phi}/sin(\widehat{\phi}+\widehat{\theta})=sin\widehat{\phi}/sin[(2\pi+(m+n)\widehat{\phi})/n]$, if we use equation 5. The first derivative of this function is as follows:\\
$A'=f'(\widehat{\phi})=\{cos(\widehat{\phi})sin[(2\pi+(m+n)\widehat{\phi})/n]-[(m+n)/n]sin(\widehat{\phi})cos[(2\pi+(m+$\
$n)\widehat{\phi})/n]\}/\{sin[(2\pi+(m+n)\widehat{\phi})/n]\}^2$\\
If we assume that $A'>0$ then this is equivalent to:\\
$g(\widehat{\phi})=tan[(2\pi+(m+n)\widehat{\phi})/n]-[(m+n)/n]tan(\widehat{\phi})>0$, of which if we take the first derivative,we have:\\ $g'(\widehat{g})=[(m+n)/n]\{sec[(2\pi+(m+n)\widehat{\phi})/n]\}^2-[(m+n)/n]\{sec(\widehat{\phi})\}^2>0$ which is true, because:\\
$|sec[(2\pi+(m+n)\widehat{\phi})/n]|>|sec(\widehat{\phi})|$ for $\widehat{\theta}+\widehat{\phi}<\pi-\widehat{\phi}$ or $A<1$. So all together we have $A'>0$ for $0<A<1$ or for $0<\widehat{\phi}<(n-2)\pi/(m+2n)$ from equation 5 and $\widehat{\theta}+\widehat{\phi}<\pi-\widehat{\phi}$, also $A\rightarrow 0$ when $\widehat{\phi}\rightarrow 0$, from equation 8 we have $\lambda \neq 0$ for $0<\kappa<I$ and the graph of equation 8 is always concave up for $0<\kappa<I$, therefore for small enough values of  $\widehat{\phi}$ there will always be pairs of values of $\kappa$ and $\lambda$ satisfying equations 6 and 8, such as in figure 3c. More specifically, for a given value of $\widehat{\phi}$, two pairs of $\kappa$ and $\lambda$ exist (two solutions), apart from the case where the curve of equation 6 is tangent to the curve of equation 8 giving the maximum value of $\widehat{\phi}$ and its related $\kappa$ value (double solution), as in figure 3c1. This is the case of the maximum value of $\widehat{\phi}$ which can be obtained for a given pair of $n$ and $m$ values and in this figure we have $m=2$, $n=12$, maximum $\widehat{\phi}=41.2degr.$ and $\kappa=0.232$. In Table 4 maximum values of $\widehat{\phi}$ and their related $\kappa$ values for specific values of $n$ and $m$ are presented and this table can be easily extended with the appropriate use of a graphic tool such as Desmos. It is useful to note that there are no closed triangle spiral cases for $0<n<3$, because of equation 5 and of the condition $0<\widehat{\theta}<\pi$.\\

\textbf{Part 3.} From Table 2 we have that  $H=sin\widehat{\phi}/(sin(\widehat{\phi}-\widehat{\theta})+sin\widehat{\theta})$, which is equation 17 for $k=1$. If we consider $\widehat{\phi}$ as a variable and the rest of the entities as constants with given value, then $H=f(\widehat{\phi})$, where $\widehat{\phi}<\pi$. The first derivative of the above function is as follows:\\

$f'(\widehat{\phi})=sin\widehat{\theta}(cos\widehat{\phi}-1)/(sin(\widehat{\phi}-\widehat{\theta})+sin\widehat{\theta})^2<0$ and $0<\widehat{\phi}<\pi$\\

For $\widehat{\phi}=0$ and following Hopital's rule, we have the following:\\

$lim(f(\widehat{\phi}))=[sin\widehat{\phi}]'/[(sin(\widehat{\phi}-\widehat{\theta})+sin\widehat{\theta})]'=1/cos\widehat{\theta}>1$ where $\widehat{\phi}\rightarrow0$,\\

Additionally for $\widehat{\phi}=\widehat{\theta}$ we have $H=f(\widehat{\phi})=f(\widehat{\theta})=1$. From all the above we can deduce that:

$H>1$ when $0<\widehat{\phi}<\widehat{\theta}$ and $H<1$ when $\widehat{\theta}<\widehat{\phi}$\\

\newpage

\textbf{Appendix 3}

\begin{center}
\textbf{TABLE 3}\
\end{center}
\begin{center}

\begin{tabular}{|c|r|r|r|r|r|r|r|r|r|r|}
  \hline
  figure  & \(n    \) & \(m        \)  & \(\widehat{\phi}    \) & \(\widehat{\theta}\) & \(\kappa\) & \(\lambda \) & \(I_{\lambda} \)& \(I_{\kappa} \) & \(\widehat{d}  \) \\
  \hline
  2b & 13 & 1 & 14.60 & 28.80 & 0.484 & 0.946 & 12 & 1 & 331.2  \\ \cline{2-10}
  2c & 13 & 2 & 15.00 & 30.00 & 0.486 & 0.895 &  6 & 1 & 165    \\ \cline{2-10}
  2d & 11 & 1 & 24.00 & 30.54 & 0.512 & 0.941 & 10 & 1 & 329.46 \\ \cline{2-10}
  3a & 12 & 2 & 20.00 & 33.33 & 0.820 & 0.967 &  - & - &   -    \\ \cline{2-10}
  3e &  5 & 1 & 21.00 & 67.80 & 0.718 & 0.936 &  4 & 1 & 292.2  \\ \cline{2-10}
  4a & 10 & 1 & 20.00 & 38.00 & 0.94  & 0.994 &  9 & 1 & 322    \\ \cline{2-10}
  4d & 10 & 9 & 20.00 & 18.00 & 0.94  & 0.945 &  1 & 1 &  38    \\ \cline{2-10}
  5  &  8 &13 & 20.00 & 12.50 & 0.78  & 0.668 &  3 & 5 & 137.5  \\ \cline{2-10}
  \hline
\end{tabular}
\end{center}
$ $\\

\begin{center}
\textbf{TABLE 4 [$\widehat{\phi}$ maximum (degrees)-$\kappa$]}\
\end{center}
\begin{center}

\begin{tabular}{|c|r|r|r|r|r|r|r|r|r|r|}
  \hline
  m-n  & \(3    \) & \(4        \)  & \(5    \) & \(6 \) & \(7 \) & \(8 \) & \(9 \) \\
  \hline
  1  & 15.8-0.187 & 25.8-0.172 & 32.6-0.157 & 37.5-0.151 & 41.3-0.149 & 44.2-0.142 & 46.5-0.134 \\ \cline{2-8}
  2  & 11.1-0.328 & 18.8-0.299 & 24.4-0.277 & 28.6-0.264 & 31.8-0.256 & 34.4-0.254 & 36.6-0.245 \\ \cline{2-8}
  3  &  8.5-0.441 & 14.8-0.390 & 19.5-0.365 & 23.2-0.346 & 26.2-0.339 & 28.5-0.329 & 30.6-0.324 \\ \cline{2-8}
  4  &  7.0-0.507 & 12.3-0.468 & 16.4-0.445 & 19.6-0.419 & 22.3-0.402 & 24.5-0.395 & 26.4-0.380 \\ \cline{2-8}
  5  &  5.9-0.571 & 10.5-0.524 & 14.1-0.493 & 17.0-0.471 & 19.5-0.457 & 21.5-0.445 & 23.2-0.437 \\ \cline{2-8}
  6  &  5.1-0.613 &  9.1-0.570 & 12.4-0.539 & 15.1-0.517 & 17.3-0.505 & 19.1-0.490 & 20.8-0.480 \\ \cline{2-8}
  7  &  4.5-0.651 &  8.1-0.609 & 11.1-0.578 & 13.5-0.553 & 15.6-0.540 & 17.3-0.534 & 18.8-0.518 \\ \cline{2-8}
  8  &  4.0-0.681 &  7.3-0.642 & 10.0-0.611 & 12.2-0.588 & 14.2-0.569 & 15.8-0.558 & 17.2-0.550 \\ \cline{2-8}
  \hline
\end{tabular}
\end{center}

$ $\\

\textbf{Appendix 4}

Let equation 5 to be related with the old system and the equivalent equation of the new system (of which the branches contra-rotate) is as follows (according to equation 7 and having as coefficients the unknown variables $X$ and $Y$ instead of $n$ and $m$):\\

$X\widehat{\overline{\theta}}=X(\widehat{\theta}-\widehat{\phi})=2\pi-Y\widehat{\phi}$\\

Since the parameter $X$ is the only one to be multiplied by $\widehat{\theta}$, the equivalence of the two equations (equation 5 and the above) dictates that $X=n$. From these two equations we get that $(X-Y)\widehat{\phi}=m\widehat{\phi}$, therefore $Y=n-m$.\\

\textbf{Appendix 5}

\textbf{Part 1.} Equations 36, 37 and 39 give us when co-rotating:\\
$\widehat{d}=(n\widehat{d}+\widehat{\phi})\widehat{\theta}/(2\pi)-(m\widehat{d}+\widehat{\theta})\widehat{\phi}/(2\pi) \Rightarrow
n\widehat{\theta}-m\widehat{\phi}=2\pi$\\
which is equation 5, therefore equation 39 is correct. Similarly we obtain equation 40 from equations 37 and 38 (branches contra-rotating).\\

\textbf{Part 2.} Equation 37 gives us (branches co-rotating):\\

$\widehat{d}=(I_{\kappa} 2\pi - \widehat{\theta})/m$ and from equation 36 we get $I_{\lambda}=(n((I_{\kappa} 2\pi - \widehat{\theta})/m)+\widehat{\phi})/(2\pi)$, therefore
$I_{\lambda}=(n/m)I_{\kappa} - (n\widehat{\theta}-m\widehat{\phi})/(2m\pi)$, which together with equation 5 gives us $n I_{\kappa} - m I_{\lambda} =1$ (equation 41 with sign at the right hand side). Similarly from equations 37, 38 and 7 we get equation 41 (contra-rotating). In the case of figure 5 (branches contra-rotate), the equations 37 and 38 give us the angle $\widehat{d}=\widehat{A_{0,3}SA_{5,0}}=137.5degr.$ by rotating clockwise.\\

If we rotate in the anticlockwise (opposite) direction, we get $\widehat{d'}=\widehat{A_{0,3}SA_{5,0}}=360-137.5=222.5degr.$, the equivalent 37 and 38 equations become: $2\pi I'_{\kappa} + \widehat{\theta}=m\widehat{d'}$ and $2\pi I'_{\lambda} - \widehat{\phi}=n\widehat{d'}$, since $I'_{\kappa}=m-I_{\kappa}$ and $I'_{\lambda}=n-I_{\lambda}$, so these give us $\widehat{d'}=I'_{\lambda}\widehat{\theta}+I'_{\kappa}\widehat{\phi}$ as equation 40 and $nI'_\kappa-mI'_\lambda=-1$ as equation 41 with the - sign at the right hand side, related uniquely with the versions of equations 37 and 38 mentioned above.\\
In co-rotating branches we have equations 36 and 37 : $2\pi I'_\kappa+\widehat{\theta}=m \widehat{d'}$ and $2\pi I'_\lambda+\widehat{\phi}=n\widehat{d'}$, so $nI'_\kappa-mI'_\lambda=-1$ ($\widehat{d'}$ and $\widehat{d}$ in opposite directions).\\

\textbf{Part 3.} From equation 42 we get: $r=(1/\lambda)^{I_\lambda} \kappa^{I_\kappa}$ which, when combined with equation 5 gives us: $r=k^{-(m/n)I_\lambda}\kappa^{I_\kappa}=\kappa^{(nI_\kappa-mI_\lambda)/n}$ and because of equations 41 and 5, this becomes $r=\kappa^{1/n}=\lambda^{1/m}<1$.\\

\textbf{Part 4.} For spirals with co-rotating branches, if equations 6 and 8 are combined by equating their $\lambda$ expressions, the following equation can be obtained :

\begin{eqnarray*}
\kappa^{(m+n)/n}sin\widehat{\theta}- \kappa^{m/n}sin(\widehat{\theta}+\widehat{\phi})+sin\widehat{\phi}=0
\end{eqnarray*}

The above equation (taking in account that $\lambda=\kappa^{1/n}$), is equivalent to the equation 10, in [2, p14].

Also for spirals with contra-rotating branches, similarly to the above, we have:

\begin{eqnarray*}
\kappa^{m/n}sin\widehat{\phi}+\kappa sin\widehat{\theta}-sin(\widehat{\theta}+\widehat{\phi})=0
\end{eqnarray*}

From equations 37 and 38, we get : $sin(n\widehat{d})=sin\widehat{\phi}$, $sin(m\widehat{d})=-sin\widehat{\theta}$ and $sin((m-n)\widehat{d})=-sin(\widehat{\phi+\theta})$ which together with $r=\kappa^{1/n}$, make the above equation equivalent to the equation 5 at [2, p.7].\\

\end{document}